\documentclass[%
 reprint,aip,
superscriptaddress,
 amsmath,amssymb,
]{revtex4-2}
\usepackage{graphicx}
\usepackage{dcolumn}
\usepackage{amsmath,amssymb,epsfig}\usepackage{commath}
\usepackage{subcaption}\usepackage{lineno,float}\usepackage{changes}
\newcommand{\beq}{\begin{equation}}
	\newcommand{\eeq}{\end{equation}}

\usepackage{xcolor}

\usepackage{enumerate}
\newcommand{\btab}[2]{\begin{tabular}{#1}#2\end{tabular}}
\def\sm{\small}\def\rb{\raisebox}
\def\ig{\includegraphics}
\def\bmip{\begin{minipage}{\textwidth}}\def\emip{\end{minipage}}
\def\huga#1{\begin{gather} #1 \end{gather}}

\def\hual#1{\begin{align} #1 \end{align}}

\def\eex{\hfill\mbox{$\rfloor$}}

\newtheorem{theorem}{Theorem}[section]
\newtheorem{lemma}[theorem]{Lemma}
\newtheorem{remark}[theorem]{Remark}
\newtheorem{example}[theorem]{Example}
\newtheorem{exercise}[theorem]{{}\hs{-1mm}}
\def\bexe{\begin{exercise}}\def\eexe{\eex\end{exercise}}
\def\bsol{\begin{solution}}\def\esol{\eex\end{solution}}
\def\bexa{\begin{example}}\def\eexa{\end{example}}
\def\brem{\begin{remark}}\def\erem{\end{remark}}
\def\bthm{\begin{theorem}}\def\ethm{\end{theorem}}
\def\blem{\begin{lemma}}\def\elem{\end{lemma}}
\def\bcor{\begin{corollary}}\def\ecor{\end{corollary}}
\def\bdefi{\begin{definition}}\def\edefi{\end{definition}}
\newcommand{\R}{{\mathbb R}}
\newcommand{\C}{{\mathbb C}}
\newcommand{\Z}{{\mathbb Z}}

\def\pdep{{\tt pde2path}}
\def\bsub{\begin{subequations}}\def\esub{\end{subequations}}

\newcommand{\bsmm}{\begin{large}\left(\begin{smallmatrix}}
		\newcommand{\esmm}{\end{smallmatrix}\right)\end{large}}

\def\CD{{\cal D}}  

\def\CO{{\cal O}}

\def\del{\delta}
\def\om{\omega}
\def\noi{\noindent}\def\ds{\displaystyle}\def\al{\alpha}
\def\pa{{\partial}}\def\lam{\lambda}
\def\Om{\Omega}
\def\dd{\, {\rm d}}\def\ri{{\rm i}}\def\er{{\rm e}}
\newcommand{\bi}{\begin{itemize}}\newcommand{\ei}{\end{itemize}}
\newcommand{\ben}{\begin{enumerate}}\newcommand{\een}{\end{enumerate}}
\newcommand{\bci}{\begin{compactitem}}\newcommand{\eci}{\end{compactitem}}
\newcommand{\bcen}{\begin{compactenum}}\newcommand{\ecen}{\end{compactenum}}
\newcommand{\bce}{\begin{center}}\newcommand{\ece}{\end{center}}
\newcommand{\reff}[1]{(\ref{#1})}

\newcommand{\hs}[1]{{\hspace{#1}}}\newcommand{\vs}[1]{{\vspace{#1}}}

\def\eps{\varepsilon}

\newcommand{\barr}{\begin{array}}\newcommand{\earr}{\end{array}}
\newcommand{\bpm}{\begin{pmatrix}}\newcommand{\epm}{\end{pmatrix}}
\newcommand{\bsm}{\left(\begin{smallmatrix}}
	\newcommand{\esm}{\end{smallmatrix}\right)}
\newcommand{\ba}{\begin{array}}\newcommand{\ea}{\end{array}}

\def\p2p{{\tt pde2path}}
\def\ind{{\rm ind}}
\def\cc{\rm c.c.}
\def\hot{\rm h.o.t.}

\usepackage{hyperref}\usepackage{cleveref}
\renewcommand{\today}{14 January 2024}

\makeatother
\begin{document}

\preprint{AIP/123-QED}

\title[]{Time-dependent localized patterns in a predator-prey model}

\author{Fahad Al Saadi}
\email{fahad.alsaadi@mtc.edu.om}
\affiliation{Department of Systems Engineering, Military Technological College, Muscat, Oman}
\affiliation{Department of Engineering Mathematics, University of Bristol, Bristol, UK}
\author{Edgar Knobloch}%
\affiliation{Department of Physics, University of California, Berkeley, CA 94720, USA}%
\author{Mark Nelson}%
\affiliation{School of Mathematics and Applied Statistics, University of Wollongong,
	Wollongong, NSW 2522, Australia}%
\author{Hannes Uecker}%
\affiliation{Institut f\"ur Mathematik, Universit\"at Oldenburg, D26111 Oldenburg, Germany}%

\date{\today}

\begin{abstract}
	Numerical continuation is used to compute solution branches in a two-component
	reaction-diffusion model of Leslie--Gower type. 
        Two regimes are studied in detail. In the first, the homogeneous
	state loses stability to supercritical spatially uniform oscillations, followed
	by a subcritical steady state bifurcation of Turing type.
	The latter leads to spatially localized states embedded in an oscillating background
	that bifurcate from snaking branches of	localized steady states. Using
	two-parameter continuation we uncover a novel mechanism whereby disconnected segments
	of oscillatory states zip up into a continuous snaking branch of time-periodic localized
	states, some of which are stable. In the second, the homogeneous state loses stability
	to supercritical Turing patterns, but steady spatially localized states embedded either
	in the homogeneous state or in a small amplitude Turing state are nevertheless present.
	We show that such behavior is possible when sideband Turing states are strongly
	subcritical and explain why this is so in the present model. In both cases the
	observed behavior differs significantly from that expected on the basis of a 
	supercritical primary bifurcation.
\end{abstract}

\maketitle

\begin{quotation}
Coupled reaction-diffusion equations describe a multitude of physical processes
ranging from morphogenesis to intracellular dynamics, catalysis and models of
vegetation cover in dryland ecosystems. In two-species models a spatially
uniform state may lose stability to a time-independent Turing pattern or
a spatially uniform oscillation, depending on parameters. We study the interaction
of these two instabilities when a primary supercritical oscillatory instability 
is followed in close succession by a subcritical Turing bifurcation. We show
that the spatially localized states associated with the latter inherit Hopf
bifurcations from the uniform state leading to localized states embedded in an 
oscillating background and show that states of this type can be {\em zipped up}
by varying a second parameter into a continuous snaking branch of such time-periodic
states, thereby demonstrating that these states may exhibit behavior analogous to
that of time-independent localized states. We also compute a 
snaking branch of steady localized states in a situation where the primary 
Turing bifurcation is supercritical, and explain this unexpected behavior 
in terms of subcritical sidebands.
\end{quotation}

\section{Introduction}
The term homoclinic snaking (HS), coined in Ref.\cite{CW99}, refers to branches of steady
localized solutions (LS) of a spatially extended pattern-forming system that oscillate
back and forth across an interval in parameter space. These oscillations reflect
the growth of the structure, with each oscillation responsible for the addition
of one wavelength of the pattern on either side of the LS. Thus the snaking
interval corresponds to the presence of a multitude of coexisting LS of
different lengths. The basic mechanism behind HS was explained by Pomeau\cite{pomeau} 
and relies on the pinning of the fronts at either end of the LS to the pattern within
it,\cite{knobloch15} a process that can be analyzed by beyond-all-orders 
asymptotics as described in Refs.\cite{chapk09,dean11,deWitt19} and references therein.
In general one finds two intertwined snaking branches, one corresponding to reflection-symmetric
states that peak on the symmetry axis and the other corresponding to states that dip on
the symmetry axis. As shown in Ref.\cite{beck} LS may also be organized in a stack of disconnected isolas
instead of the two continuous intertwined solution branches. Localized states are generally only
found when the primary Turing bifurcation is subcritical, resulting in bistability between the trivial
state and the pattern state, leading to the interpretation of LS
as segments of the coexisting Turing pattern embedded in the competing homogeneous background. 
It is expected
that stable LS are only present when the competing states are both stable. Exceptions are 
known in the form of slanted snaking, for instance in the presence of global coupling\cite{FCS07}
or a neutral large scale mode,\cite{dawes08,thiele,beaume} for which localized states are found
even in the supercritical case, with no bistability. See Ref.\cite{thiele} for a more precise
thermodynamic interpretation of this behavior and Refs.\cite{uwsnak14,w18,UW18a, KUW19} for HS
of ``patterns within patterns'', which are related to more complicated bistabilities.

The present work concerns two fundamental questions: \\
(i) Do time-dependent LS exhibit similar behavior to that of the well-studied steady states?\\
(ii) Do LS exist in systems with local coupling and no conserved quantity when the primary 
bifurcation is supercritical?

Time-dependent LS are of two basic types, traveling pulses and localized standing oscillations. The
former are known to exhibit snake-like behavior in appropriate parameter regimes and may be organized
in a stack of isolas consisting of 1,2,... equispaced traveling pulses as found in natural doubly diffusive
convection,\cite{lojacono} or in a single snaking branch as found for a three-species reaction-diffusion
(RD) system in Ref.\cite{KUY21} For localized standing oscillations the situation is much less clear, although
these are expected to behave in a similar fashion as the steady LS on account of their shared reflection
symmetry. In this work we present a simple two-species RD system which resolves both issues. Specifically,
we show that this system contains time-periodic snaking states obtained by ``zipping up'' tertiary
branches of localized oscillatory states, with the snaking structure of steady LS serving as a template or
backbone. Moreover, we show that in a different parameter regime the same system exhibits steady LS even
when the primary Turing bifurcation is supercritical and explain why.

The zipping up process is of particular interest. We show that it is a consequence
of the nonlinear interaction between a supercritical Hopf bifurcation of the trivial
state and a nearby subcritical Turing instability, and in particular of the presence of
tertiary Hopf bifurcations on both the Turing branches and the associated LS
inherited from the primary Hopf bifurcations of the trivial state. As such this
behavior appears characteristic of the interaction between the Turing and Hopf
instabilities\cite{dewit96,MWBS97,TMV09,TBMV11,TMBMV13,p2pDMV} and we demonstrate
its presence in two different RD systems.

The first system we study is a Leslie-Gower prey--predator model taking the form 
\cite{Zhou14,Zhou19} 
\begin{subequations}\label{eq1}
	\begin{alignat}{2}
		\pa_t u & = D \pa_x^2 u{+}f(u,v)
		& \equiv &D \pa_x^2u{+}u\bigl(a{-}u{-}bvh(u,v)\bigr),  \label{1a11} \\
		\pa_t v&=\pa_x^2 v+g(u,v)
		& \equiv & \phantom{\,\eps}\pa_x^2 v +\del v \bigl(1 - vh(u,v)\bigr),   \label{1b111} 
	\end{alignat}
\end{subequations}
where $u$ and $v$ are the prey and predator densities, respectively. 
The system is posed on a 1D domain $\Om=(-\ell,\ell)$ with large $\ell$ and Neumann boundary
conditions (NBC) $\pa_x u|_{\pa\Om}=\pa_x v|_{\pa\Om}=0$ on both components. 
We choose 
$$
\mbox{$h(u,v)=\ds 1/\left( (\alpha u+1)( \beta v+1)\right)$}, 
$$
known as Bazykin's functional response.\cite{bazykin1998}
The parameters $a$, $b$, $\alpha$, $\beta$ and $\delta$ are all positive, 
and the prey diffusion constant $0<D\ll 1$.  The parameter 
$\del$ plays the role of a time scale for the predator ODE, and hence can 
be used to induce Hopf instabilities, as discussed in, e.g., \cite{CPVB23} 
for a modified version of Bazykin's system. 

The system \reff{eq1} has four homogeneous steady states, 
namely $s_1=(0,0)$, $s_2=(0,\left( 1-\beta \right)^{-1})$, $s_3=(a,0)$, 
and the coexistence state 
\huga{\label{s4}
	\text{$s_4=(u^*,v^*)$: 
		$u^* \equiv a-b$, $\ds v^*\equiv\frac{\alpha u^* +1}{1-\beta -\alpha \beta u^*}$,}
}
present provided $a>b$ and $1-\beta > \alpha \beta u^*$. 
The steady states of \reff{eq1}, in particular $s_4$, and their 
instabilities have been analyzed in detail in Ref.\cite{Zhou19}, and 
spatio-temporal solutions have been obtained via direct numerical 
simulation. 

Of particular interest are the Hopf and Turing bifurcations from $s_4$ as 
parameters vary, in particular the parameter $b$.  We use the 
toolbox \pdep\ \cite{p2pbook} to numerically continue the primary
solution branches bifurcating from $s_4$ together with the secondary 
and tertiary branches
that bifurcate from them, resulting in a large multiplicity of steady 
and time-periodic LS.
For simplicity, our study is restricted 
to two parameter ranges. In (i), present at low $b$ in Fig.~\ref{bifur1}, 
we study the interaction
of LS with Hopf instabilities of the homogeneous background $s_4$. 
In (ii), present at higher $b$ in Fig.~\ref{bifur1}, we find 
a novel example of HS where (large amplitude) LS embedded in $s_4$
transition to LS embedded in a small amplitude background pattern
arising in a supercritical Turing bifurcation from $s_4$. 
This state is similar to those in Ref.\cite{KUW19}, where bistability between 
small and large amplitude patterns of the same wavelength 
was created artificially by considering a cubic-quintic-septic Swift-Hohenberg equation.
Similar two-scale structures have been reported in rotating convection,\cite{beaume}
rotating Couette flow\cite{salewski19} and in natural binary fluid convection.\cite{tumelty23}

In both our regimes, (i) and (ii), a key feature is that it is {\em not the primary} 
Turing bifurcation that determines the existence of stable localized 
patterns. Instead, periodic patterns which bifurcate further 
away from the primary Turing bifurcation become ``most stable'', and penetrate
farthest into the (subcritical) 
parameter range of stable homogeneous steady states. Thus, beyond 
the result that rather simple and natural two-component RD systems 
can generate more complicated HS than models such as the 
(quadratic or cubic) SH equation, one important lesson
is that it may be necessary to go beyond the 
first few Turing branches to understand the possible 
multitude of stable steady states, both spatially extended and spatially localized.

The paper is organized as follows. In Sec.~2 we summarize the linear stability properties of
the trivial state $s_4$. This is followed in Sec.~3.A by a detailed study of case (i), and in
Sec.~3.B by case (ii). In Sec.~4 we show that analogous behavior occurs in the Gilad-Meron
model of dryland vegetation, and use this result to argue that the zipping up mechanism
uncovered here is a generic process. Brief conclusions follow in Sec.~5.

\section{Linear stability analysis}
\begin{figure*}[ht]
\btab{lll}{{\sm (a)}&\qquad{\sm (b)}&{\sm (c)}\\
\hs{-3mm}\ig[width=0.48\textwidth]{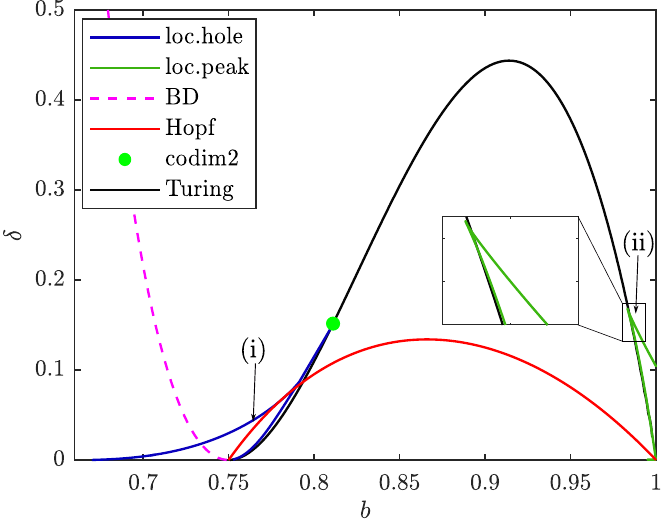}&
\qquad\rb{32mm}{\btab{l}{\ig[width=42mm,height=32mm]{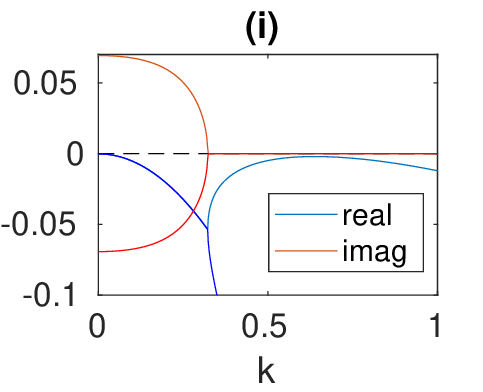}\\
\ig[width=42mm,height=32mm]{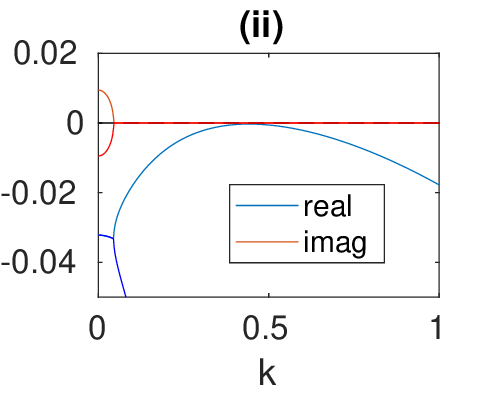}}}&
\rb{30mm}{\btab{l}{\ig[width=42mm,height=33mm]{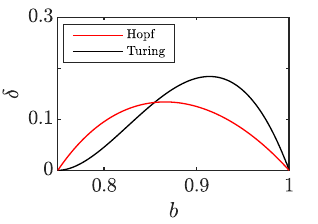}\\
\ig[width=42mm,height=30mm]{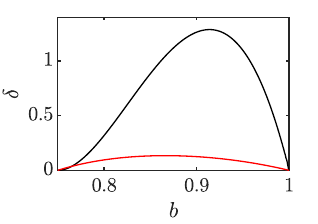}}}
}

\vs{-2mm}
\caption{{(a) Hopf and Turing instability curves for the homogeneous solution
    $s_4$ of \eqref{eq1} in the 
    $(b,\delta)$ plane when $(a,\alpha,\beta,D)=(1,2,0.5,0.029)$ together
    with the existence regions for localized patterns. (b) Spectra $\lambda(k)$
    at the locations indicated in (a). (c) Bifurcation curves for $D=0.07$
    (top) and $D=0.01$ (bottom) showing the dependence of the Turing
    instability (black curves) on the parameter $D$. The Hopf bifurcation
    (red curve) is independent of $D$. 
\label{bifur1}}}
\end{figure*}

The linear stability of $s_4$ is described by
\begin{equation}\label{eq:linear}
(u,v)^T_t =\left[J|_{s_{4}}+\CD\pa_{x}^2\right](u,v)^T, 
\end{equation}
where $(u,v)^T$ denotes the {\em perturbation} of $(u^*,v^*)^T$, $J|_{s_4}$ is the reaction Jacobian 
at $s_4$ and $\CD\equiv\mbox{diag}(D,1)$ is the diffusion matrix. The Fourier ansatz 
$(u,v)^T=(A,B)^T e^{\ri k x+\lambda t}+\cc$, where $\cc$ denotes the 
complex conjugate of the preceding terms, leads to a dispersion 
relation $\det(J-k^2\CD-\lambda(k)I)=0$
connecting the growth rate $\lambda(k)$ to the assumed wave number $k$.
Thus, $s_4$ is (linearly) 
stable if $\Re\lambda(k) <0$ for all $k$, while $\Re\lam(k)=0$ 
at some {\em critical wave number} $k_c\in\R$ indicates, subject to transversality conditions,
the onset of instability and hence a bifurcation. 
These are classified as longwave ($k_c=0, \Im\lam(0)=0$), 
Turing ($k_c\ne 0$, $\Im\lam(k_c)=0$), or Hopf 
($k_c=0, \Im\lam(0)\ne 0$). A fourth possibility is the  
so-called wave instability ($k_c\ne 0$, $\Im\lam(k_c)\ne 0$), 
but this cannot occur in two-component RD systems. 
The Turing instability arises at $\delta=\delta_T$, where 
\hual{
  \bigl[&\alpha\left( {\alpha}\beta\delta_T D{-}2 \right) {u^*}^{2}{+}
 \left(  \left( 2\beta{-}1 \right) \alpha \delta_T D{+}a\alpha{-}1
 \right) u^*\notag\\
 &{+} \left( \beta{-}1 \right) \delta_T D\bigr]^2
=
4u^* \left( u^*\alpha{+}1 \right)^2 \delta_T D \left( 1{-}\beta \left( u^*\alpha{+}
1 \right)  \right),\label{Turineq}
}
provided the critical wave number $k_T$ given by
$\ds k_T^2=\frac {\alpha \left( \alpha\beta\delta_T D{-}2 \right) {u^*}^
{2}{+}\left( \left( 2\beta{-}1 \right) \alpha \delta_T D+a
\alpha{-}1 \right) u^*{+} \left( \beta{-}1 \right)\delta_T D}{ 2\left(
  u^*\alpha+1 \right) D}$
is real, while the Hopf instability arises at $\delta=\delta_H$, where
\begin{equation}\label{Hopfcon}
   \delta_H= {\frac { \left( a-b \right)  \left( 1+ \alpha \left( a-2\,b \right) 
 \right) }{ \left( \alpha\beta\, \left( a-b \right) +\beta-1 \right) 
 \left( \alpha\, \left( a-b \right) +1 \right) }}, 
\end{equation}
independently of the diffusion coefficient $D$, provided the Hopf frequency 
$$\omega_H= \sqrt{-\delta_H (a-b)[\alpha \beta (a-b) +\beta -1]}$$
is real. Figure \ref{bifur1}(a) shows these instability curves in 
the $(b,\delta)$ parameter plane for 
$(a,\alpha,\beta,D)=(1,2,0.5,0.029)$, while (b) shows 
the eigenvalue curves $k\mapsto \lam_{1,2}(k)$ 
at the two points labeled (i) and (ii), which will also be our 
starting points for the numerical bifurcation analysis. Panel (c) shows
the dependence of the Turing and Hopf curves on the parameter $b$ for two
different values of $D$.

\brem\label{krem}{\rm a) The Fourier ansatz 
$(u,v)^T=(A,B)^T e^{\ri k x+\lambda t}+\cc$ 
with 
wave number $k\in\R$ applies to infinite domains $x\in\R$. 
For the numerics we have to choose a finite 
domain $x\in(-\ell,\ell)$, and we choose homogeneous 
NBCs for both $u$ and $v$. This restricts the wave number $k$ 
to $k\in \frac \pi {2\ell}\Z$. Nevertheless, 
for large $\ell$, 
this gives a rather 
dense sampling of the dispersion relation $k\mapsto \lam(k)$. 

b) In the supercritical range (i.e., after crossing the 
Turing or Hopf line), on the infinite line, we have a 
band of unstable wave numbers, and bifurcations of $2\pi/k$ 
periodic patterns for $k$ arbitrary close to $k_c$. These ``sideband 
patterns'' are unstable at bifurcation but stabilize 
at small but finite amplitude via Eckhaus bifurcations. 
This behavior is inherited for large but finite $\ell$, where 
the discrete wave numbers $k\in \frac \pi {2\ell}\Z$ 
lead to a close succession of sideband Turing, resp.~Hopf, bifurcations 
after the primary Turing, resp.~Hopf, bifurcation. 

c) The statements a) and b) hold for general Turing and Hopf instabilities. 
However, a pronounced feature of \reff{eq1} is that 
Turing branches with nearby but different $k$ may behave quite 
differently, i.e., their direction of branching is very sensitive
to $k$. As a consequence, to understand LS in \reff{eq1} (as opposed to, e.g.,
the cubic Swift--Hohenberg equation) 
it is {\em not} sufficient to determine whether the primary ($k=k_c$) bifurcation 
is sub- or supercritical, and the sidebands must be taken 
into account (see Sec.~3.2).
}\eex 
\erem 

In Fig.~\ref{bifur1}(a) we mark a codimension-two point (green dot), 
determined from
weakly nonlinear theory, where the Turing bifurcation changes from 
sub-- to supercritical. 
Moreover, we show two curves labeled ``localized holes'' (blue, 
on the left) and two labeled ``localized peaks'' (green, on the right)
delimiting the regions of existence of LS.%
\footnote{These curves were computed by fold-continuation 
of {\em selected} folds of LS; as explained below, it is in general not 
clear which LS (i.e., of which width and wave number) extends farthest in 
the $(b,\delta)$ plane, and hence these lines provide an approximate characterization only.} 
The LS in the former region
take the form of finite arrays of large amplitude downward spikes, while the
latter consist of upward spikes. 
Finally, we also show the Belyakov-Devaney (BD) transition curve (dashed purple line), 
defined by a pair of real spatial eigenvalues of the spatial 
dynamics problem linearised about the homogeneous equilibrium, each of double multiplicity.
This curve is also given by \reff{Turineq} but the corresponding wave number is now imaginary.
One major significance of the BD transition is that (standard) snaking of 
LS typically turns into ``foliated snaking'' of 
isolated spikes\cite{ponedel16,yochelis21,FahadWoods,CHAMPNEYS1998158,ALSAADI2021125014} 
when crossing the 
BD line; see Ref.\cite{VERSCH2021} for a general analysis. 
Here, this region of foliated snaking 
(to the left of the BD line and below the loc.holes line) is rather 
small and will not be studied. 
In Fig.~\ref{bifur1}(c) we show the instability curves for $D=0.01$ 
and $D=0.07$. For smaller $D$
the Turing curve expands, while for larger $D$ it shrinks (in \eqref{Turineq} $D$ only appears in the combination $\delta D$) while the Hopf curve is unaffected.
In both cases we lose the interaction of the Hopf and Turing modes on the left (near (i)), 
motivating our choice of an intermediate value $D=0.029$ as a starting point for our study.

\section{Two case studies of localized patterns}
\label{nsec}
We use numerical continuation with \pdep\ \cite{p2pbook} to explore 
patterns near the two selected instability points from 
Fig.~\ref{bifur1}. First, we consider the vicinity of 
the left point (i), where the primary instability is of Hopf type. 
However, subcritical Turing instabilities are nearby, leading 
to localized steady patterns, and we shall see that these inherit,
in some sense, the Hopf instability of the homogeneous 
background. Second, we explore the neighborhood of the right point (ii). 
Here the primary instability is a supercritical Turing instability, but
strongly subcritical sideband instabilities nearby
lead to peculiar stable localized patterns, for which the 
associated snake interpolates between homoclinics to 
the homogeneous state and homoclinics to small amplitude periodic 
patterns. For our numerical continuation we choose $\Om=(-\ell,\ell)$ with
$\ell=50$, and plot 
\hual{\label{n1}
	\|u\|_2:=\sqrt{\frac 1 {2\ell}\int_{-\ell}^{\ell} |u(x)-u^*|^2\dd x}, \quad\text{resp.}\quad\\ 
	\|u\|_2{:=}\sqrt{\frac 1{2T\ell} \int_0^T \int_{-\ell}^{\ell}
		|u(x,t){-}u^*|^2\dd x\dd t}\label{n2}
}
for steady, resp.~time-periodic states, i.e., we only use the first component
of the {\it perturbation} from $s_4$, and normalize its norm by the domain size
(and the period $T$ for time-periodic orbits). The projection $\|u\|_2$ as defined
is {\em not} a norm in the $(u,v)$ space. 

\subsection{Interaction between LS and background Hopf instability}\label{nsec1}
\begin{figure*}
\btab{l}{
\btab{ll}{{\sm (a)}&{\sm (b)}\\[-4mm]
\hs{-2mm}\ig[width=66mm]{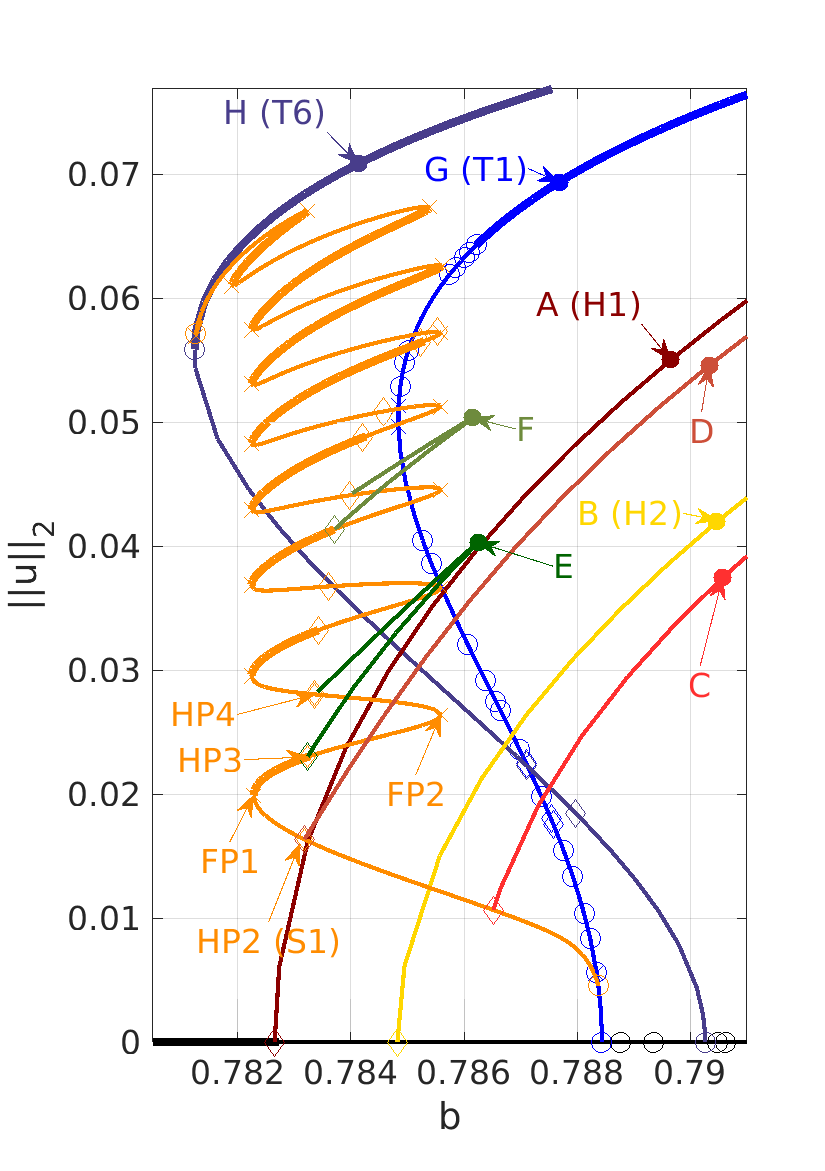}&
\hs{-2mm}\rb{50mm}{\btab{l}{\hs{10mm}\vs{0mm}
\ig[width=14mm]{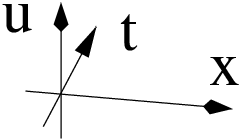}\hs{40mm}\ig[width=14mm]{xtc3}\\[-1mm]
\ig[width=53mm]{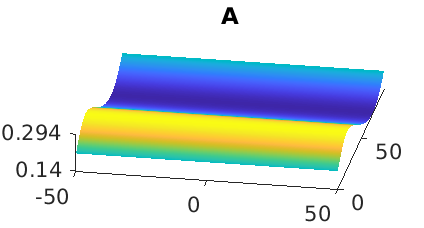}
\ig[width=53mm]{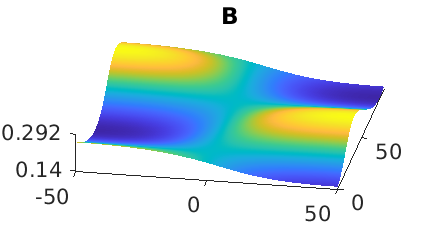}\\
\ig[width=53mm]{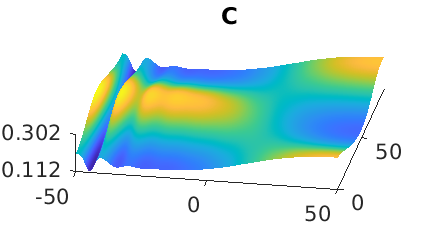}
\ig[width=53mm]{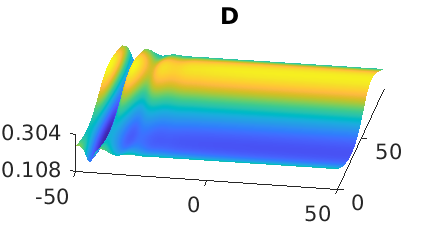}\\
\ig[width=53mm]{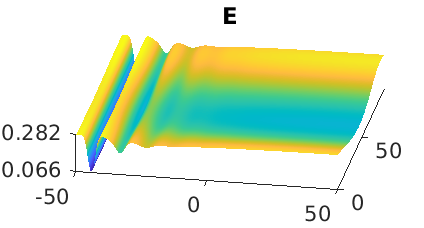}
\ig[width=53mm]{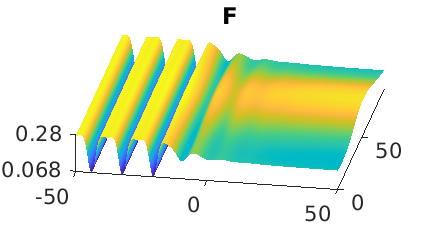}
}}}\\[-0mm]
{\sm (c)}\\[-4mm]
\hs{10mm}\ig[width=11mm]{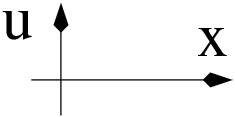}\\
\ig[width=42mm]{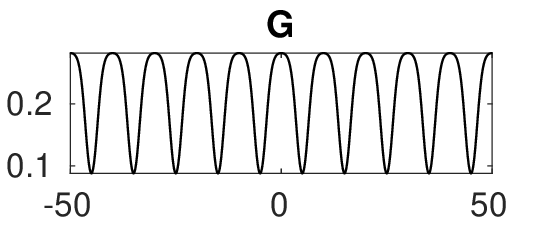}
\ig[width=42mm]{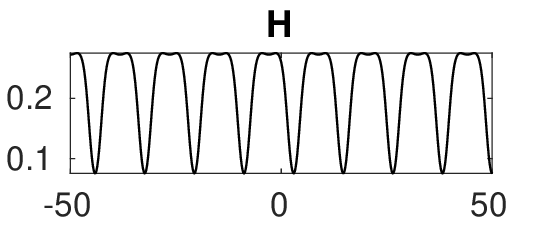}\ig[width=42mm]{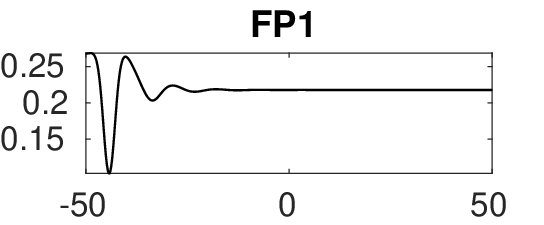}
\ig[width=42mm]{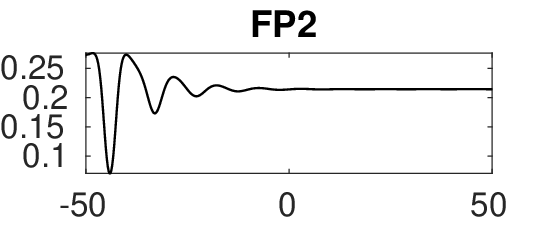}
}
\caption{(a) Bifurcation diagram for $\ell=50$ with $(\del,\alpha,\beta,a,D)=(0.07,2,0.5,1,0.029)$ showing $\|u\|_2$ defined in \reff{n1} and \reff{n2} as a function the parameter $b$ starting from $b=0.78$. (b) Sample time-periodic solutions ($u$ component only) shown in a space-time plot over one temporal period corresponding to the locations labeled in (a). (c) Sample steady solutions ($u$ component only). 
In (a) Hopf bifurcation points (HPs) are marked by $\diamond$, steady state
BPs by $\circ$, and fold points (FPs) by $\times$. Not all 
HPs, BPs and FPs are shown. For steady states, thick (thin) lines indicate stable (unstable)
solution branches. The trivial state $s_4$ corresponds to $u=0$.
Hopf branches H1 and H2 from this state are shown in brown and yellow, 
with sample solutions A and B. The primary Turing branch T1 is shown in blue;
the 6th Turing branch (T6, dark blue) reaches farthest to the left. A snaking branch of LS
(S1, orange) bifurcates from the first BP on T1 and reconnects to T6, sample solutions in (c).  
Additional HPs on S1 lead to branches with sample solutions (C-F). 
\label{f1}}
\end{figure*}

Figure \ref{f1}(a), with $\ell=50$, shows that if $b$ is increased at fixed $\del=0.07$
near location (i), then $s_4$ loses stability to a Hopf mode with $k_c=0$, i.e., to a spatially uniform
oscillation (H1, see the space-time plot at location A 
in panel (b)). This is followed
by a Hopf bifurcation to a spatially nonuniform Hopf mode with wave number $\pi/\ell$
(H2, shown in a space-time plot at location B).
Both bifurcations are supercritical, implying that the $k=0$ branch H1 is 
stable at onset, while the $k=\pi/\ell$ oscillations H2 are unstable at onset, with Floquet index 
$\ind=1$ (number of Floquet multipliers $\mu$ with $|\mu|>1$).

The next bifurcation point (BP) from $s_4$ yields the first 
Turing bifurcation (T1, bluebranch, $k=\pi/5\approx 0.628$, profile G). 
This bifurcation is subcritical, like the next 10 Turing bifurcations, 
and the ``most subcritical'' Turing branch is the 6th 
(T6, dark blue branch, $k=0.534$, profile H).%
The branch T1 undergoes a secondary
bifurcation at small amplitude, i.e., close to the primary bifurcation,
to a pair of spatially modulated states which turn into a pair of intertwined snaking branches of 
localized Turing patterns only one of which (S1, orange branch) is shown in the figure, see profiles
at successive fold points FP1 and FP2 in panel (c); 
owing to NBC these profiles may be reflected in $x=-50$ to generate a localized state on
a domain of length 200. As $\ell$ increases this secondary bifurcation moves to smaller amplitude
and collides with the primary bifurcation in the 
limit $\ell\to\infty$.\cite{bukno2007}

Along the orange branch S1 of steady LS the solution adds a wavelength near every
other fold until the available domain is filled and the snaking branch terminates on T6.
This is a consequence of the fact that the wavelength within the LS is not set by $k_c$ but
is instead set by nonlinearity and hence the parameter $b$.\cite{bukno2007,bergeon08}
Note that S1 snakes outside of the region of bistability between $s_4$
and T1, despite its origin in T1. Instead the relevant bistability region is the region
between $s_4$ and T6. On T1 and T6 (and similarly 
on all other Turing branches, not shown) there are many further BPs, 
leading to similar secondary bifurcations to localized patterns 
and snaking branches like the orange branch S1; for instance, the second BP 
on T1 yields genuine LS on $(-\ell,\ell)$, i.e., double pulse homoclinic orbits
on the domain of length 200.

The above behavior is largely consistent with the standard snaking of LS
observed in the Swift-Hohenberg equation. The present system is not variational,
however, and so time-dependent LS may be expected.\cite{burkedawes} Moreover,
in the present system snaking occurs in a region of bistability between
steady LS and smaller amplitude uniform oscillations on the Hopf branch H1.
In this region steady LS cannot be stable over unbounded domains 
since the background state $s_4$ is unstable to oscillations, but stable LS
embedded in a time-periodic background may be expected. Evidently, the
stability of various states on S1 indicated in Fig.~\ref{f1}(a) is a finite
size effect. 

Next we consider tertiary Hopf points (HP) on S1. The first two are inherited
from the oscillatory instability of $s_4$ and generate states that resemble a
``superposition'' of an LS on S1 and the Hopf modes H2 and H1, see profiles 
C and D, respectively. In the following
we refer to these states as mixed modes (MM). These time-periodic states
grow monotonically in amplitude as $b$ increases, and both are unstable. 
As we continue upward along the S1 branch we find pairs of HPs that are connected
by further branches of MMs (green). In each case the lower HP is on a stable 
part of the LS branch while the upper one is on the next unstable
part. Altogether we have 6 pairs of such HPs, all pairwise connected
by an MM branch. Similar states were found in Ref.\cite{CPVB23},
Figs.~13(h) and 14(b), using DNS for suitable values of the predator time 
scale parameter (our $\del$) by starting from a localized perturbation of the 
analog of our $(u^*,v^*)$, and 
in the Gilad-Meron model of dryland plant ecology,
also a two-species RD system (Ref.\cite{ASP23}, Fig.~9).

We think of these green branches as forming parts of an unzipped ``snake'',
broken by segments of steady LS (see below for zipping up). Sample solutions
at locations E and F are shown in panel (b).
In all cases the green segments between the lower HP and the corresponding fold point
correspond to {\it stable} MMs. 
We have checked, but do not show, that the intertwined
second LS snaking branch (with minima at $x=-50$) undergoes essentially identical
behavior with pairs of BPs straddling the right folds and generating partially stable
MM segments closely resembling the green branches in Fig.~\ref{f1}(a). These states
are also embedded in a uniform background oscillation.

There are two finite domain effects responsible for the existence of
{\em stable} steady LS connected to the background state $s_4$ in a
parameter regime in which $s_4$ is predicted to be unstable to
oscillations: as one proceeds up the S1 branch the LS grow in extent,
thereby reducing the domain occupied by $s_4$.  Figure
\ref{bifur1}(b), panel (i), shows that as the allowed minimum
wave number is pushed to higher values the Hopf bifurcation is
suppressed. We expect therefore that the tertiary HPs move to larger
values of $b$, as observed. This does not explain, however, the
presence of stable narrow LS, such as those created at FP1. In fact
FP1 lies below the primary Hopf point for $s_4$, and $s_4$ is
therefore stable at this location. The fact that this state extends
stably past the primary Hopf point is a consequence of the fact that
once an LS is present the stability of the background state no longer
reflects the stability of this state on an unbounded domain (or the
equivalent problem with NBC): the smaller oscillation scale at the
location of the front between the LS and the background $s_4$ enhances
dissipation leading to an increase in the critical $b$ for the onset of
the Hopf instability, cf.~Refs.\cite{tobias98,batiste06}. Thus the observed
behavior is a consequence of both a reduction in the effective domain
size, but more importantly of the introduction of a strongly dissipative
(non-NBC) boundary at the location of the front separating the LS and
the oscillating background.  The resulting dissipation suppresses the
oscillation amplitude in this region, as evident from the space-time
profiles E and F in Fig.~\ref{f1}(b). We remark that once the LS occupies
approximately half the domain the role of the LS and the oscillating
state changes: it is now more natural to think of a uniform
oscillation embedded in a background of a steady periodic pattern.

\subsubsection{Zipping up the TH snake}\label{zsec}

\begin{figure*}
\btab{ll}{{\sm (a)}&{\sm (b)}\\[-0mm]
\hs{-4mm}\rb{30mm}{\btab{l}{
\ig[width=40mm,height=30mm]{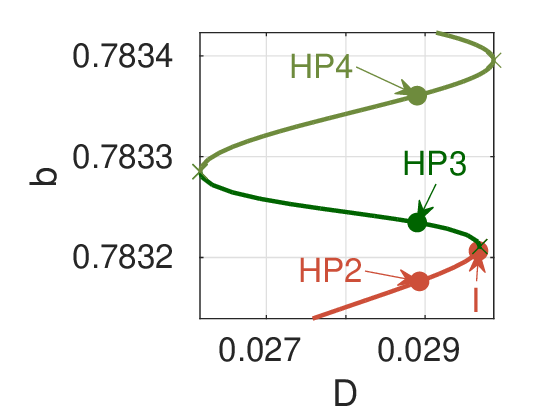}\\
\ig[width=40mm,height=30mm]{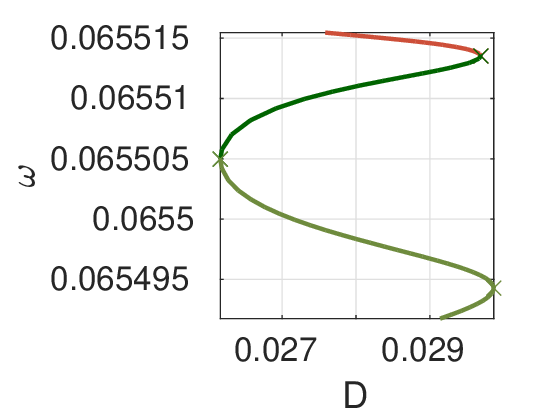}}}&
\hs{-4mm}\ig[width=58mm,height=65mm]{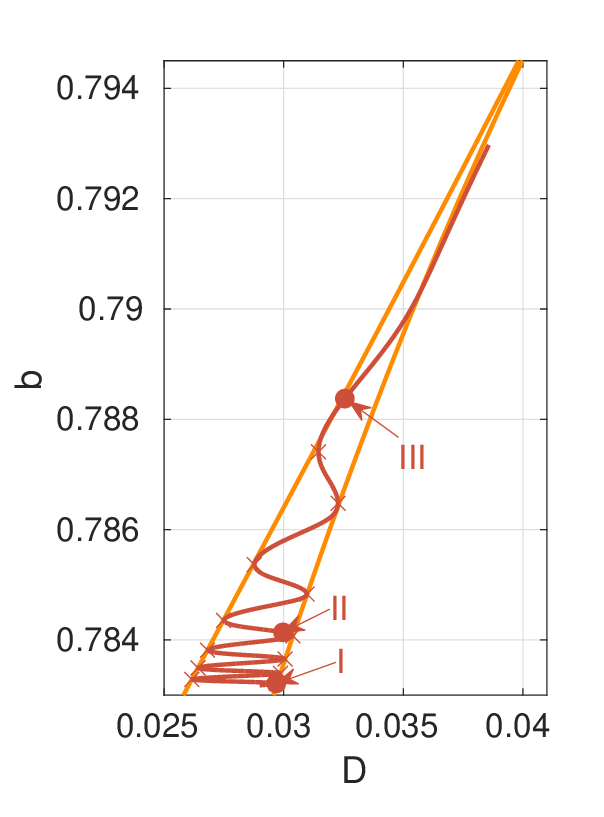}
\hs{-0mm}\rb{32mm}{\btab{l}{
\hs{6mm}\ig[width=12mm]{xc2}\\
\hs{-2mm}\ig[width=34mm,height=18mm]{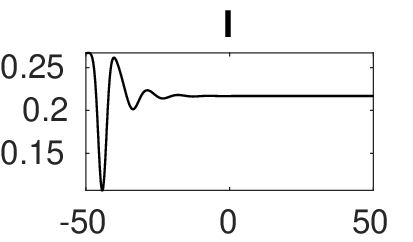}\\
\ig[width=31mm,height=18mm]{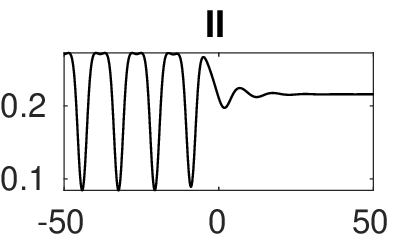}\\
\ig[width=31mm,height=18mm]{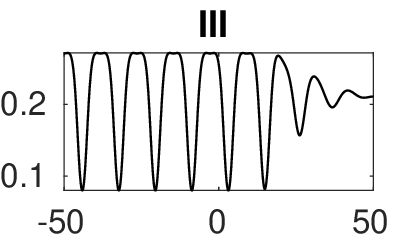}
}}}\\
\btab{ll}{{\sm (c)}&{\sm (d)}\\[-0mm]
\hs{-4mm}\ig[width=58mm]{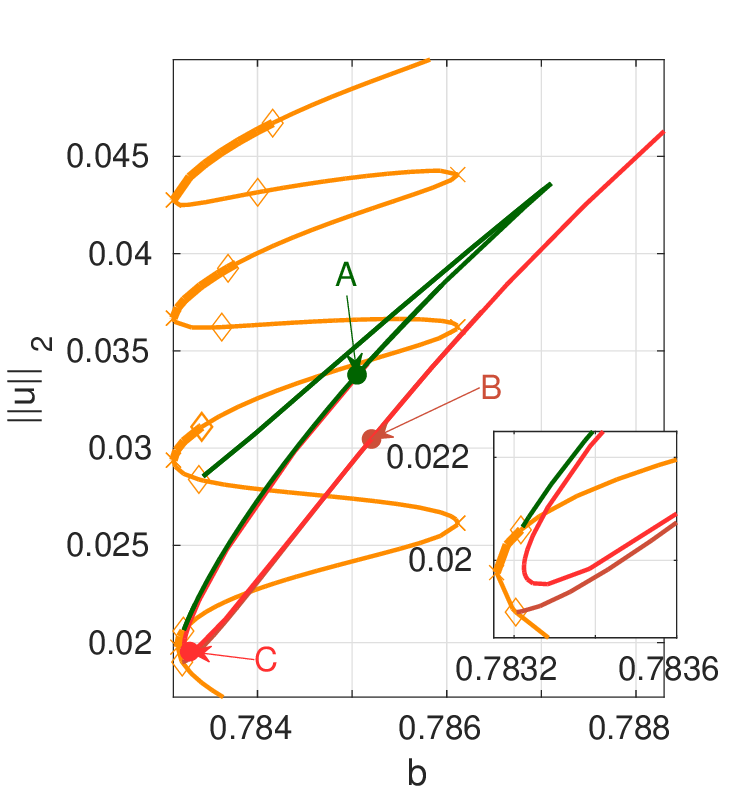}
\hs{-1mm}\rb{28mm}{\btab{l}{
\hs{10mm}\vs{0mm}\ig[width=13mm]{xtc3}\\[-1mm]
\ig[width=37mm]{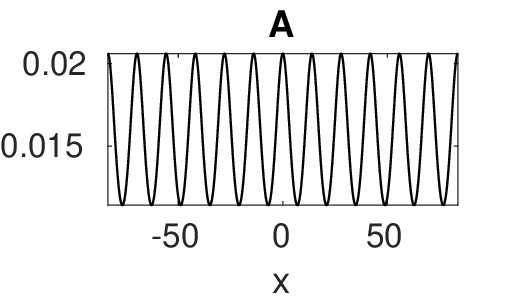}\\
\ig[width=37mm]{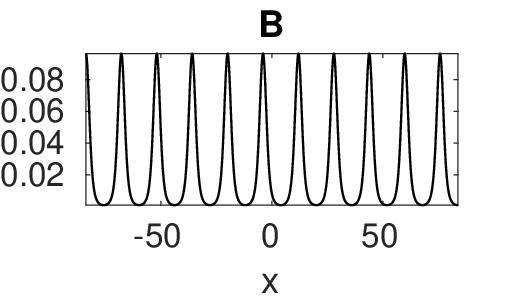}\\
\ig[width=37mm]{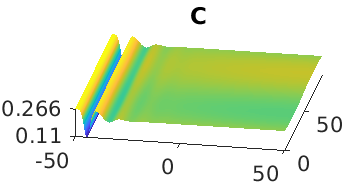}}}
&\hs{0mm}\ig[width=48mm,height=65mm]{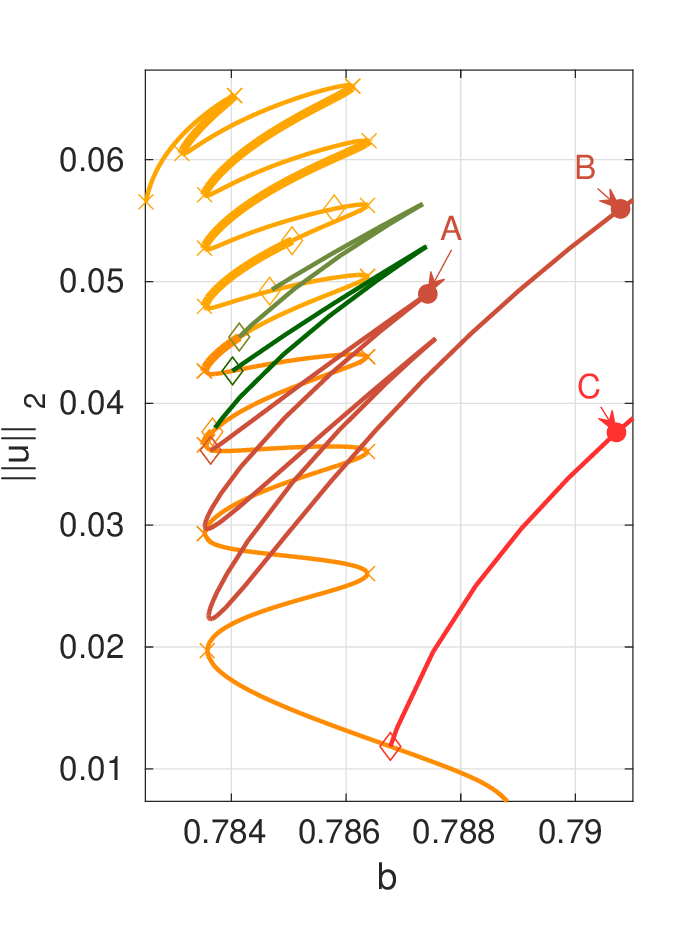}
\hs{-2mm}\rb{31mm}{\btab{l}{
\hs{10mm}\vs{0mm}\ig[width=13mm]{xtc3}\\[-1mm]
\ig[width=37mm]{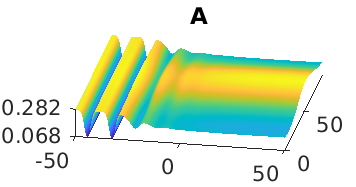}\\
\ig[width=37mm]{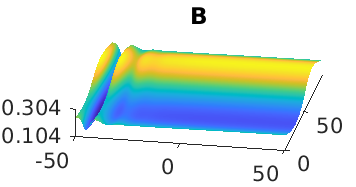}\\
\ig[width=37mm]{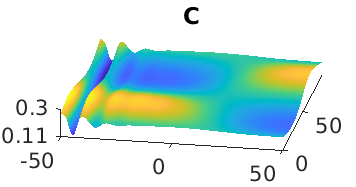}}}}
\vs{-2mm}
\caption{{\sm Zipping up the Turing-Hopf segments from Fig.~\ref{f1} into a 
Turing-Hopf snake on increasing $D$. (a) HPC of HP2 (brown), HP3 (dark green) 
and HP4 (light green) from Fig.~\ref{f1} showing the HP location in the $(D,b)$
plane, together with the corresponding Hopf frequency $\om$. (b) HPC of HP2 (brown),
and FPC of FP2 and FP3 (orange) together with sample solutions I-III used to initialize
continuation in $b$ in panels (c) and (d). (c) Bifurcation diagram
for $D=0.02965$ corresponding to point I in panel (b); the inset shows the  
LS branch (orange) near the first fold, shortly before HP2 and HP3 annihilate
(at $D=D_1\approx 0.02968$) together with the disconnected MM branches (brown and green). 
The red branch (inset) shows the MM branches at $D=0.0297$, i.e., after reconnection. As
a result the large $b$ MM branch now connects to the point HP4 on S1. To generate the
red branch we take point A from the green branch at $D=0.02965$, set $D=0.0297$, and
continue in $b$. (d) Bifurcation diagram for $D\approx 0.03$ corresponding to point II
in panel (b). Stability on S1 is indicated by thick lines; stability on the other
branches is not indicated (see text).
\label{f3b}}}
\end{figure*}

\begin{figure}[ht]
\bce 
\hs{-4mm}\ig[height=72mm]{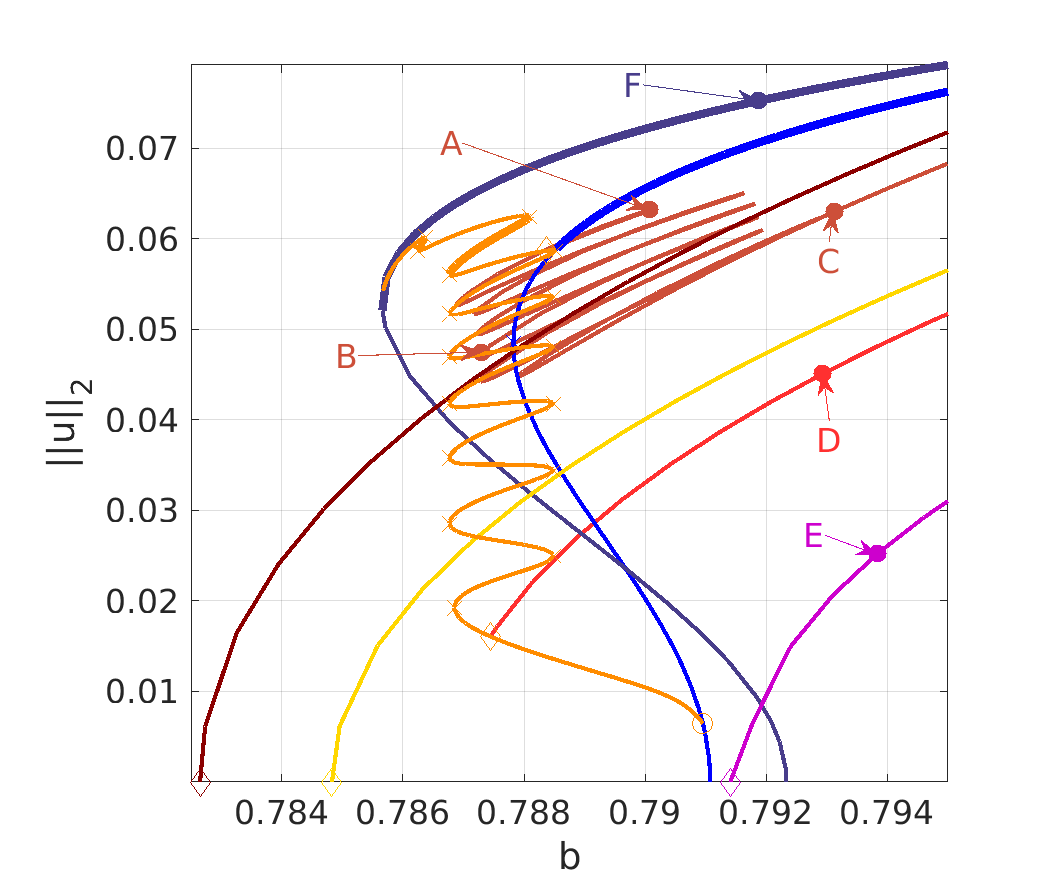}
\hs{-2mm}\rb{33mm}{\btab{l}{
\hs{10mm}\vs{0mm}\ig[width=12mm]{xtc3}\\[-1mm]
\ig[width=39mm]{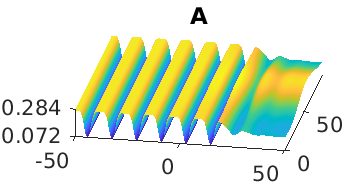}\\
\ig[width=39mm]{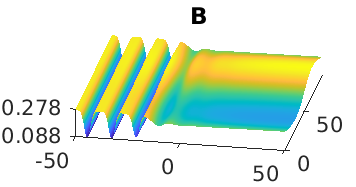}\\
\ig[width=39mm]{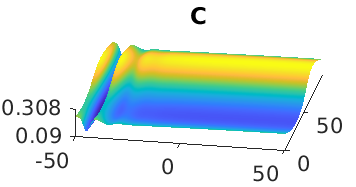}}}
\hs{-2mm}\rb{33mm}{\btab{l}{
\ig[width=39mm]{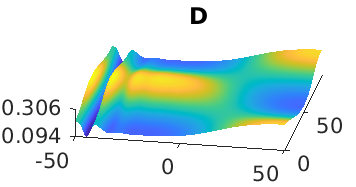}\\
\ig[width=39mm]{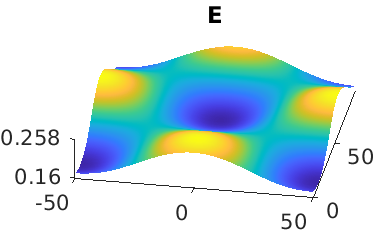}\\
\ig[width=37mm]{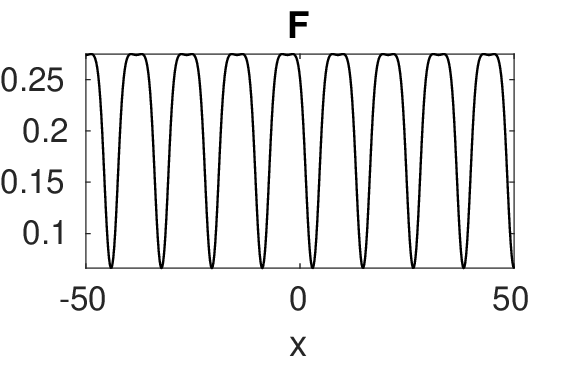}}}
\ece

\vs{-4mm}
\caption{{\sm As for Figs.~\ref{f3b}(c,d) but for $D=0.033$, corresponding to point III in Fig.~\ref{f3b}(b): 
    the whole Turing-Hopf snake has zipped up, generating a snaking branch of time-periodic states
    terminating near the top of S1 (solutions A-C). Stability is not indicated. The only other HP left
    on S1 is the one at the bottom leading to the (red, unstable) MM branch with a nonuniform background
    oscillation (solution D here and solution C in Fig.~\ref{f3b}(d)).
    The figure also shows other branches at $D=0.033$, namely H1 (brown), H2 (yellow), H3 (violet, with sample
    solution E), as well as T1 (blue) and T4 (4th Turing branch from $s_4$, dark blue with sample solution F),
    the termination of S1 for this value of $D$. 
\label{f3b2}}}
\end{figure}

Figures \ref{f3b} and \ref{f3b2} show how the (green) snake segments
reconnect or ``zip up'' into a full snaking branch of time-periodic
states by varying the parameter $D$.  In Fig.~\ref{f3b}(a), obtained
from Hopf-point continuation (HPC) of HP2, HP3 and HP4 on the S1
branch in Fig. \ref{f1}(a), we show the Hopf point positions $b$ and
the corresponding frequencies $\om$ as $D$ varies. As $D$ increases,
the location of HP2 moves up in $b$ while the location of HP3 moves
down, resulting in a collision at the fold FP1 of S1 when
$D\approx 0.02968$.  Similarly, as $D$ decreases, the location of HP3
moves up in $b$ while the location of HP4 moves down, resulting in a
collision at FP2 when $D\approx 0.026$.  The same happens for HPC of
other pairs of Hopf points such as HP5, HP6 and HP7, HP8, etc.  Put
differently, if we follow HP2 beyond the first collision we obtain the
brown curve in Fig.~\ref{f3b}(b), a {\it snake of Hopf points}. The
folds in this curve correspond to successive collisions of Hopf
points, HP2 with HP3 at FP1, HP3 with HP4 at FP2 etc.  As a result the
HP snake is trapped between, say, the FPC of the fold points FP2
(upper boundary) and FP3 (lower boundary) of Fig.~\ref{f1}(a) since the
right and left fold points of S1 are almost aligned. Sample solutions
corresponding to the locations I, II and III in (b) are shown
alongside. Importantly, with increasing $D$, the wedge containing the
S1 branch becomes exponentially narrow in $|D-D_c|$ with
$D_c\approx 0.043$, cf.~Refs.\cite{chapk09,dean11,deWitt19}, and overall
the structure leans to the right.

Figures~\ref{f3b}(c,d) describe steps in the zipping process.  In
panel (c) we recompute the orange, green and brown curves in
Fig.~\ref{f1}(a) for $D\approx 0.02965$. We see that the Hopf points
HP2 and HP3 have moved very close to FP1, but the branches bifurcating
from HP2 (brown, going to large $b$) and HP3 (green, reconnecting to
HP4 on S1) are still distinct. For $D=0.0297$, however, the points HP2
and HP3 on S1 have annihilated one another, resulting in the
reconnection of the brown and green branches and the creation of the
red branch (panel (c), inset). As a result HP2 and HP3 are absent and
the MM branch now originates from HP4 higher up the snaking S1
branch. A further increase in $D$ leads to a collision between HP4 and
HP5 at FP3 and the process repeats, resulting in the brown  
curve in Fig.~\ref{f3b}(d) with a footpoint at HP6, just prior to the next
collision at FP5.

In Fig.~\ref{f3b2}, computed for $D=0.033$, only the (former) HP12
remains, and consequently all the green
branches have zipped up and the large $b$ MM branch (brown) now
connects to a footpoint at B. Figure~\ref{f3b2} also shows the other
branches at $D=0.033$, with the same colors as in Fig.~\ref{f1}(a)
where applicable. The orange LS snake now reconnects to the 4th Turing
branch (dark blue, profile F);
a third Hopf branch from $s_4$ is now present and is shown in
violet (profile E). The details of the switching of the termination point of S1
from T6 to T4 are expected to proceed via T5 and to resemble the process
described in Ref.\cite{bergeon08} for the Swift-Hohenberg equation but have
not been studied in detail in the present case.

The MM segments on the other LS snake intertwined with the S1 snake zip up
via a very similar process (not shown). As a result in Fig.~\ref{f3b2}
there are in fact two intertwined snaking branches of time-periodic MM
states, only one of which is shown.

We refer to the brown branch from Fig.~\ref{f3b2} as a Turing-Hopf (TH)
snake. The space-time profiles show that the TH snake
recapitulates the behavior of the S1 snake, albeit in a larger and
shifted parameter interval. In particular, as one proceeds up the TH
snake, the (almost) stationary core of the solution adds a new
wavelength during every back and forth oscillation of the branch,
restricting the oscillating background to an ever smaller part of the
spatial domain. The progressive nature of the zipping up process is
a consequence of this fact which results in HP collisions that are staggered in $D$
(Fig.~\ref{f3b}(b)). This follows from the fact that as the central LS
structure of the oscillations expands it becomes harder to excite
oscillations in the rest of the available domain. In larger domains,
therefore, the zipping process may be faster. The process itself has
one major consequence: it progressively erases the stability of the
steady LS states and ``replaces'' these states by coexisting stable
segments of time-periodic states that extend over a wider parameter
interval. These states resemble the steady LS in their center but
are embedded in an oscillating background.
The light red branch that bifurcates from HP1 on S1 (solutions C in
Fig.~\ref{f3b}(d) and D in Fig.~\ref{f3b2}) is qualitatively unaffected 
by the small changes in $D$ required to zip the TH snake.

\begin{figure}
\bce
\btab{l}{{\sm (a)}\\
\hs{-7mm}\ig[width=100mm]{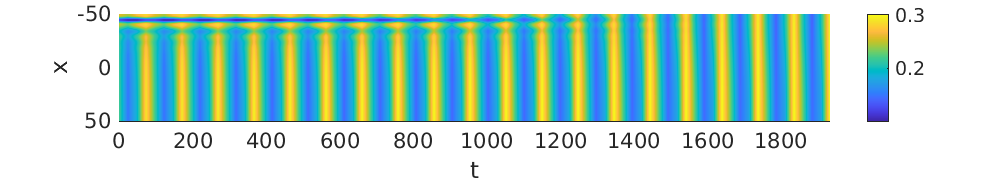}\\
{\sm (b)}\\
\hs{-7mm}\ig[width=100mm]{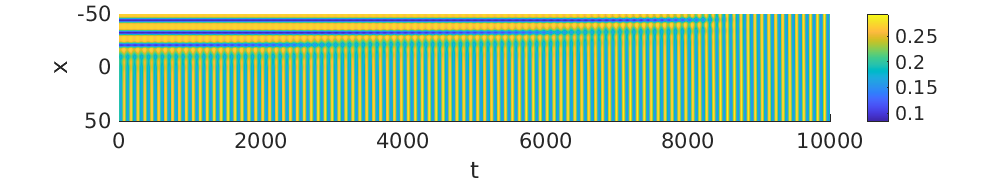}\\
{\sm (c)}\\
\hs{-7mm}\ig[width=100mm]{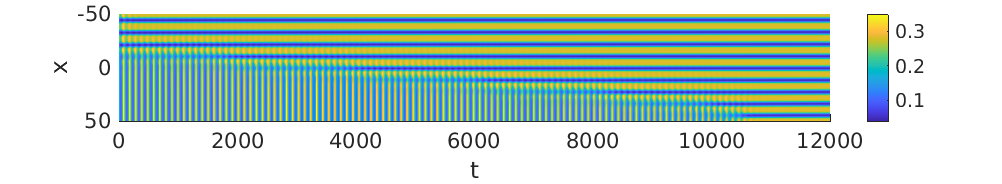}\\
{\sm (d)}\\
\hs{-7mm}\ig[width=100mm]{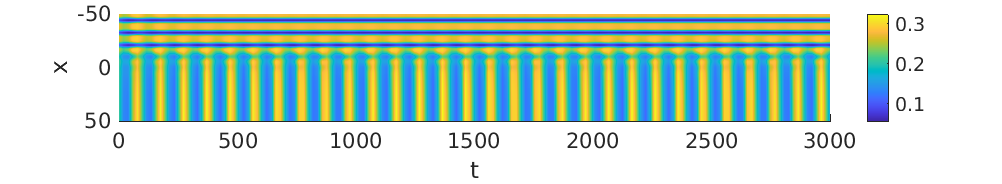}
}
\ece

\vs{-6mm}
\caption{{\sm DNS from marked points (+small perturbation) 
in Fig.~\ref{f3b2}(a). (a) Starting from C with fast 
convergence to H1. (b) Starting from B with slow convergence to H1. 
(c) Starting from B but setting $b=b+0.015$ (right of the TH snake), 
yielding (stick--slip) convergence to T4. (d) Starting from B but setting 
$b=b+0.007$, yielding convergence to a stable MM in the TH snake. 
\label{fdns}}}
\end{figure}

\brem\label{reconrem}{\rm a) The collision of HP2 with HP3 on S1
	(Figs.~\ref{f3b}(a,c)) occurs when their two frequencies and respective
	eigenfunctions become identical.  The resulting double Hopf
	bifurcation with 1:1 resonance corresponds to a nilpotent
	bifurcation, described by a normal form derived and studied in
	Ref.\cite{langford}. However, resonant bifurcations of this type are
	codimension--3 bifurcations, and are therefore not expected when
	only two parameters such as $b$ and $D$ are varied. This is because
	two parameters are required to generate a double Hopf bifurcation in
	the first place, and a third parameter is required to tune their
	frequencies into resonance. No reconnection takes place when the
	frequencies are incommensurate -- the bifurcations pass through one
	another. The required resonance forces the resonant double Hopf
	bifurcation to occur at the fold FP1, and similarly for all the
	subsequent collisions and associated reconnections. Nominally, a 1:1
	double Hopf bifurcation at FP1 is a codimension--4 event but the
	fact that this forces $d\omega_1/dD=d\omega_2/dD=\infty$
	(Fig.~\ref{f3b}(a)), in addition to $\omega_1=\omega_2$, makes the
	problem codimension--2, as documented in Fig.~\ref{f3b}(b). Simply
	put, reconnection cannot occur generically unless it takes place via
	the fold points FP on S1.
	
	b) In Fig.~\ref{fdns} we briefly look at DNS from selected points in
	and near the TH snake from Fig.~\ref{f3b}(c) obtained from initial
	conditions corresponding to these points but increasing the
	parameter $b$ by a small amount. Panel (a) shows a space-time representation of a
	solution starting from point C; the solution evolves into a stable
	H1 oscillation. Panel (b) starts from point B and also evolves into
	H1 albeit much more slowly. Panel (c) also starts from point B but
	with $b$ increased by $\Delta b=0.015$; this time the state evolves
	into a stable spatially periodic state on the T4 branch.  Finally,
	panel (d) starts from point B but with $\Delta b=0.007$; the
	solution converges to a stable MM state on a green MM branch right
	above that shown in the figure. 
	Each of these four panels shows an advancing or retreating front
	between two states (steady or oscillating) that exhibits stick-slip
	motion, much as in the Swift-Hohenberg equation.\cite{bukno2006} 
	Evidently in this region, both
	before and after the zipping up process, the system exhibits extreme
	multistability that increases with the domain size $\ell$.
	
	c) Hopf instabilities of LS can also appear in a form other than as
	``background oscillations'' in Fig.~\ref{f1}, namely in the form of
	``breathing peaks'', where the background remains at rest but the
	localized state oscillates. See, e.g., Refs.~\cite{ALSAADI23,FahadWoods}.
        However, such breathing peaks have not been found in \reff{eq1}.
	
	d) For background Hopf wave number $k_H\ne 0$ the background
	oscillations can organize either into traveling or standing waves
	depending on parameters.\cite{knobloch86}  However, since the onset
	of the $k_H\ne 0$ oscillations is preceded by the $k_H=0$ onset,
	both of these states are expected to be unstable. In a related
	problem arising in binary fluid convection driven by the Marangoni
	effect\cite{assemat08} the background is unstable with respect to
	$k_H\ne 0$ oscillations and fills with standing waves, with the
	constituent left-traveling waves dominant to the right of the LS and
	right-traveling waves dominant to the left -- a consequence of the
	effective boundary conditions at the location of the front(s)
	separating the oscillations from the nearly steady LS.

        e) Throughout this paper we have computed solutions on the half-domain
        imposing reflection symmetry in $x=-\ell$ to generate even solutions on
        a domain of length $4\ell$. Because of this procedure we can only compute
        snaking branches of even states. However, just as in the case of steady
        LS there is a family of of time-periodic states with odd symmetry in $x=-\ell$
        and hence a pair of intertwined branches of even and odd time-periodic states
        on a domain of length $4\ell$.
} \eex\erem

\subsubsection{Continuation of Hopf and fold points in $\delta$}\label{hpc2sec}

\begin{figure*}
\btab{lll}{{\sm (a)}&{\sm (b)}&{\sm (c)}\\[-0mm]
\hs{-4mm}\rb{6mm}{\btab{l}{\ig[width=54mm]{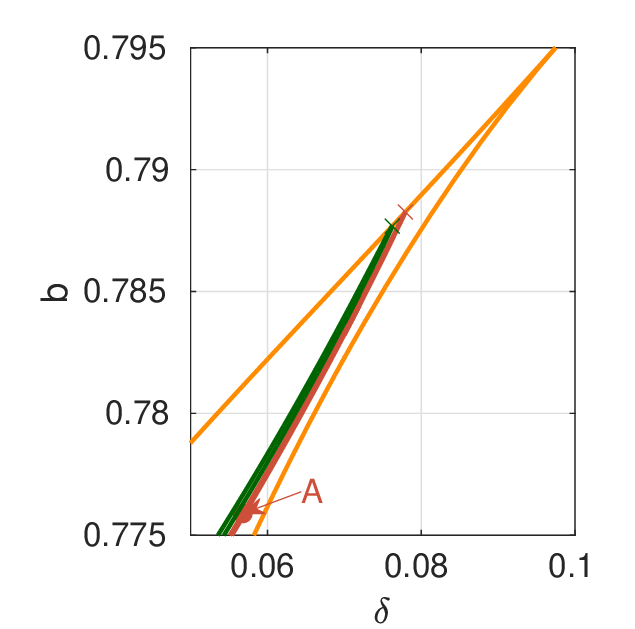}\\
\hs{10mm}\ig[width=42mm]{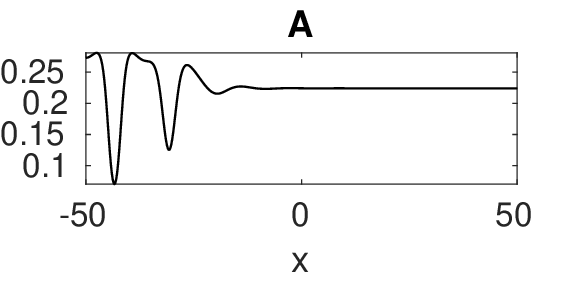}}}&
\hs{-1mm}\rb{-30mm}{\ig[width=66mm]{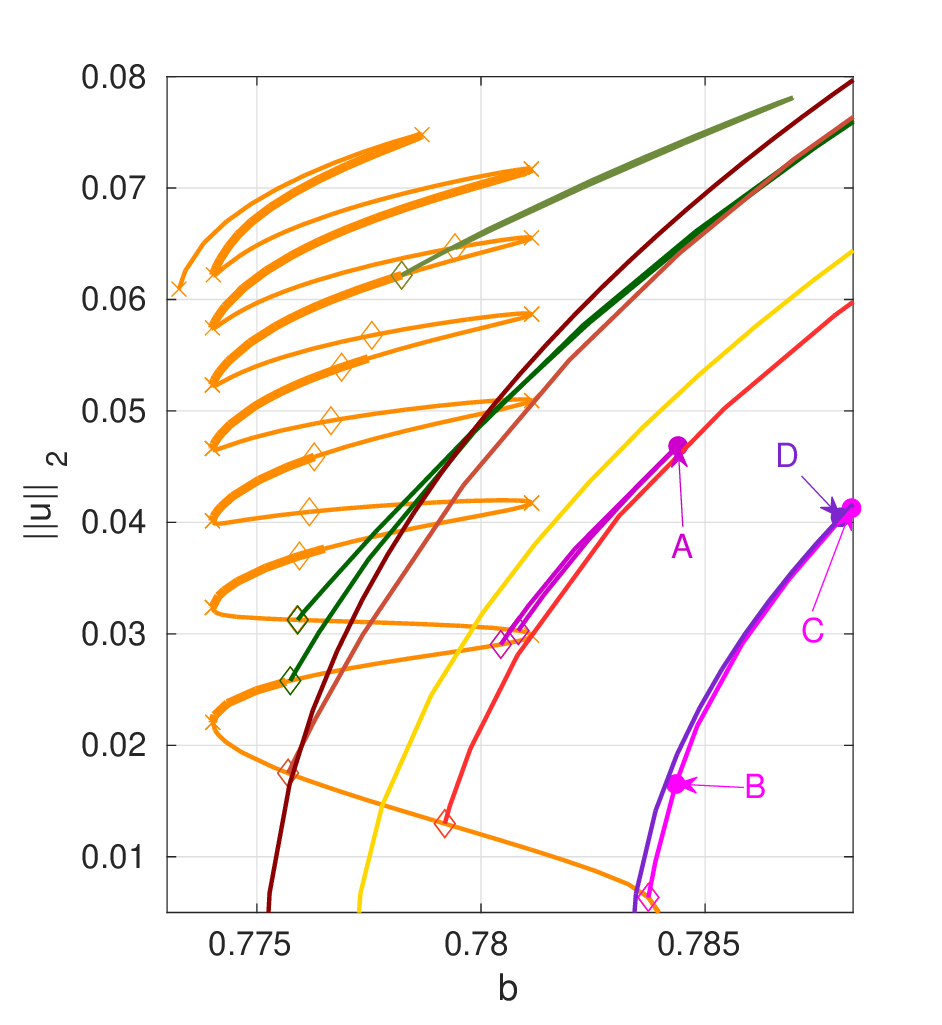}}
&\hs{-1mm}\rb{5mm}{\btab{l}{\hs{10mm}\vs{0mm}\ig[width=12mm]{./xtc3}\\[-1mm]
\ig[width=40mm]{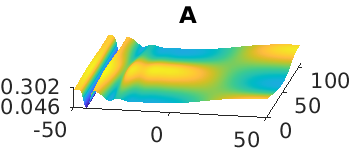}\\
\ig[width=40mm]{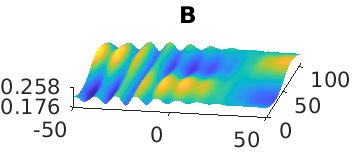}\\
\ig[width=40mm]{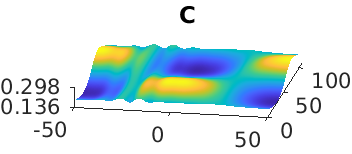}\\
\ig[width=40mm]{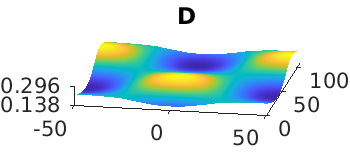}
}}}

\vs{-2mm}
\caption{{\sm As for Fig.~\ref{f1} but showing the dependence on the parameter $\del$
    instead of $D$. (a) The $(\del,b)$ plane, with FPC of FP1 and FP2 (orange), and
    HPC of HP2 (brown) and HP7 (green), with a sample solution profile at 
    $\del\approx 0.056$. (b) Bifurcation diagram at $\del=0.056$; the branches
    already present in Fig.~\ref{f1} (same colors here) behave as before but
    extend farther to the right, and there are two new (unstable) Hopf branches:
    one secondary mixed mode branch (pink, profiles B and C) with background
    wavenumber $2\pi/100$, and the corresponding connection on S1 (purple,
    profile A), analogous to the green branches with $k=0$ background oscillation.}
\label{f3c}}
\end{figure*}

Figure~\ref{f3c}(a) shows continuation results of bifurcation points
from Fig.~\ref{f1} as functions of the parameter $\del$ when
$(\al,\beta,a,D)=(2,0.5,1,0.029)$.  In this case the HPC never reaches
the left folds and therefore no zipping up takes place. Panel (b)
shows the continuation in $b$ from point A, $\del\approx0.056$.  There
are two main differences compared to Fig.~\ref{f1}(a): there are now
{\em three} primary Hopf branches that bifurcate from $s_4$ at ``low''
$b$, i.e., in the $b$ range of the S1 snake, with the branch H3 (violet)
of spatially dependent background oscillations with wave number
$k=2\pi/100$ (solution D in panel (c) or solution E in
Fig.~\ref{f3b2}) ``moving in'' to lower $b$. As a consequence, at the
bottom of the S1 snake there is a new HP giving rise to the pink
branch of mixed mode oscillations in the form of a superposition of the
LS from S1 and the $k=2\pi/100$ background oscillation.
At the same time two new HPs with background wave number $k=\pi/100$ appear on
the S1 snake, inherited from H2 (yellow branch) and straddling FP2, yielding a
second connecting segment (purple) of oscillations, this time with the
background wave number $k=\pi/100$.  These MM states, like those in green,
extend far outside the S1 snaking range. On larger domains, we find
additional pairs of HPs near the right folds of the S1 snake, but here
only one is present since the growth of the LS along S1 suppresses the
$k=\pi/100$ oscillations in the background. Moreover, on larger
domains the $k\ne 0$ segments again zip up when $D$ is varied (not
shown). We anticipate that for other parameter values (and/or larger
domains) the solution structure of mixed modes will incorporate states
with yet more complex background oscillations. However, these new
mixed mode oscillations are all expected to be unstable owing to the
instability of the $k\ne 0$ background oscillations with respect to
$k = 0$ oscillations.

\subsection{Supercritical Turing instability and localized
	patterns}\label{sec2}
\begin{figure*}
\bce
\btab{l}{
\btab{ll}{{\sm (a)}&{\sm (b)}\\ 
\hs{-4mm}\ig[width=60mm]{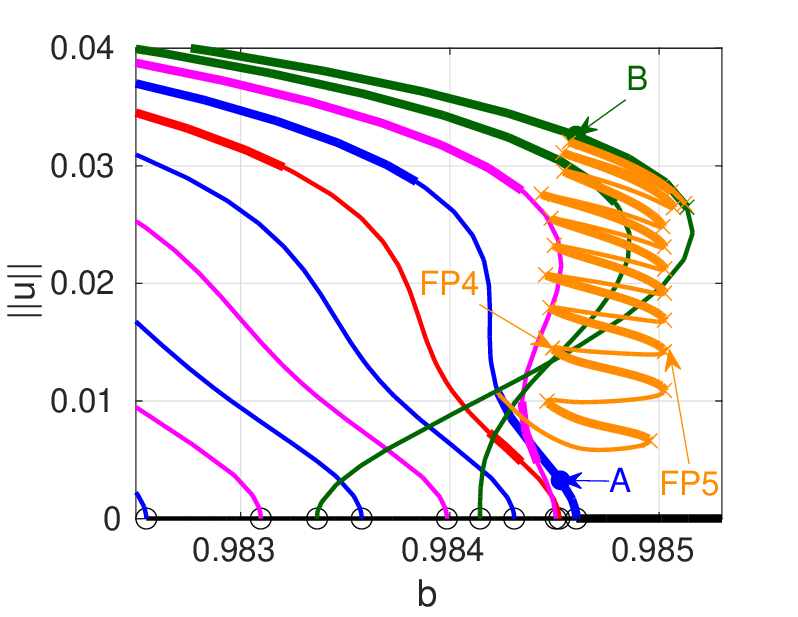}&
\ig[width=40mm]{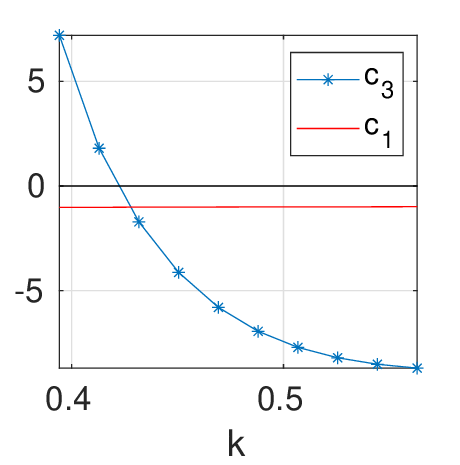}\hs{2mm}\ig[width=70mm]{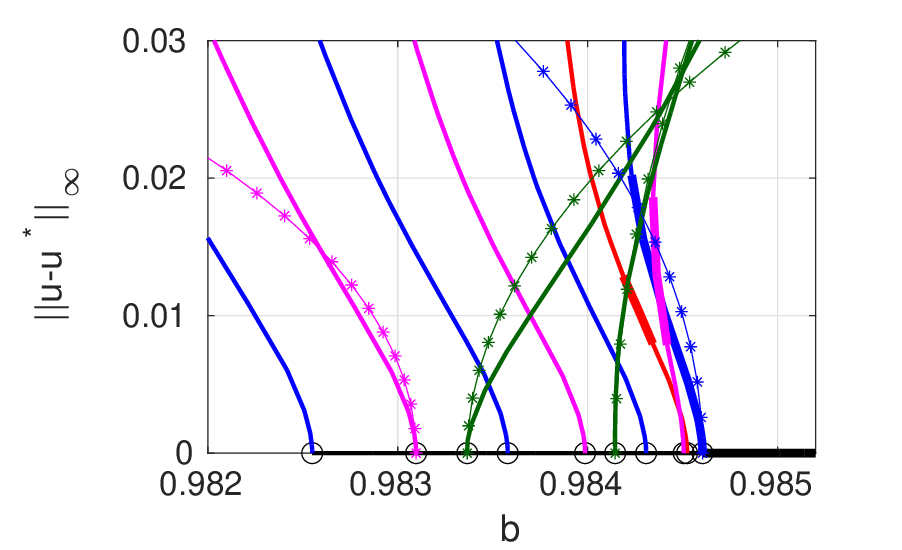}}\\[-3mm]
{\sm (c)}\\
\ig[width=42mm]{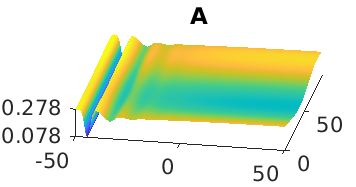}\ig[width=42mm]{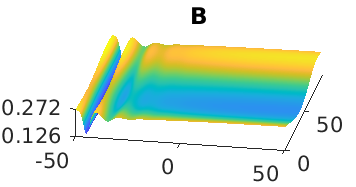}
\ig[width=42mm]{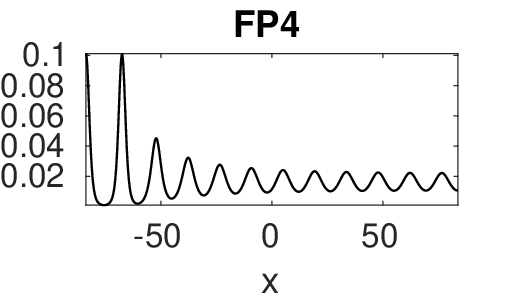}\ig[width=42mm]{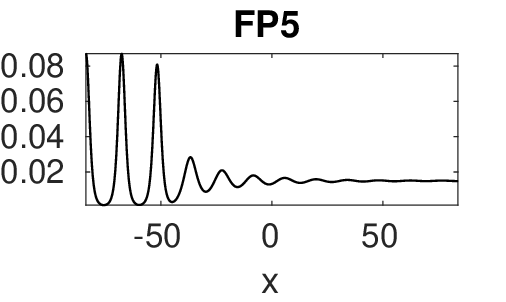}
}

\ece 
\vs{-8mm}
\caption{{\sm Bifurcation behavior in region (ii): $(\del,\alpha,\beta,a,D)=(0.16,2,0.5,1,0.029)$
    with $\ell=12\pi/k_c$ and $b$ as the primary bifurcation parameter. The onset wave number is
    $k_c=0.4505$ (T1, blue branch). (a) Bifurcation diagram showing branches T1--T9 (see Table
    \ref{Ttab} for colors) together with the S1 snake (orange).  (b) Dependence of the Landau
    coefficients on the wave number $k$ together with a comparison between the solution branches
    determined from numerical continuation (panel (a)) and the Landau approximation to selected
    branches. (c) Sample solution profiles ($u$ component only) at locations labeled in (a). }
\label{f5}}
\end{figure*}

Spatially localized patterns are usually associated with subcritical
instabilities, such as the subcritical Turing instability; see, e.g.,
Refs.\cite{Vec,yochelis21}.  However, as discussed below, LS may be found even
when the Turing bifurcation is supercritical.  In Fig.~\ref{f5} we
illustrate the behavior of \reff{eq1} on the right of the Turing
region in Fig.~\ref{bifur1}, more specifically near location (ii).
This is far from the Hopf instability of $s_4$, and indeed Hopf
instabilities no longer play a role in Fig.~\ref{f5}. Moreover, this
region is located to the right of the green point in
Fig.~\ref{bifur1}, in a region where the primary Turing bifurcation is
supercritical. In general snaking is not expected in the supercritical
regime, owing to the absence of bistability. An exception is provided
by systems with a conservation law such as the conserved
Swift-Hohenberg equation\cite{thiele} where a secondary, strongly
subcritical instability may destabilize a supercritical Turing state
at small amplitude, generating LS exhibiting (slanted) snaking in a
parameter regime with no bistability. Similar small amplitude
secondary instabilities have also been found in other conserved
systems such a rotating convection between free-slip boundaries.\cite{beaume}

Even without a conservation law, Eq.~\reff{eq1} supports LS 
in the supercritical region (ii) but for a different
reason. Figure~\ref{f5} shows that the primary Turing branch (T1,
blue) loses stability shortly after bifurcation to a (nonslanted)
branch (S1, orange) of snaking localized patterns. Near the right
folds of this structure the LS are embedded in the stable $s_4$ state
(solution FP5). However, the snake extends across the primary
bifurcation to T1 and solutions in this region are instead embedded in
a small amplitude periodic Turing pattern (solution A) as shown by
solution FP4. The transition between these two situation is smooth,
with no change in stability.  Similar states were found inrotating
plane Couette flow\cite{salewski19} and in natural binary fluid
convection\cite{tumelty23}, and were studied in Ref.\cite{KUW19} in a
model problem, the 3-5-7 Swift-Hohenberg equation; here we find that
they arise in a natural way in the RD system \reff{eq1}.

To understand the behavior shown in Fig.~\ref{f5} we computed the
Landau coefficients $c_1$ and $c_3$ that describe the small amplitude
behavior of the Turing branches. The Landau approximation for the 
branch T$_j$ with wave number $k_j$ bifurcating at $b=b_{T_j}=:b_j$
is given by 
\begin{subequations}\label{gla}
	\begin{align}
 (u,v)(t,x)&-(u,v)^*(b_j)=\eps
	A(T)\er^{\ri k_j x}\phi(k_j)\nonumber\\ & +\eps^2\biggl[\frac 1 2
	A_0(T) + A_2(T)\er^{2\ri k_j x}\biggr] +\cc+\hot, 
	\end{align}
\end{subequations}
where $\phi_j\equiv \phi(k_j)\in\R^2$ is a kernel vector of $J|_{(u,v)^*(b_j)}-k_j^2\CD$ 
(in the following normalized to $\phi_{j,1}=1$),  
$0<\eps\ll 1$ is a formal perturbation parameter, and 
$\hot$ denotes higher order terms. 
In general the complex amplitude function $A$ depends on the slow time
scale $T=\eps^2 t$ while $A_0(T)\in\R^2$ and $A_2(T)\in\C^2$ are the
modes excited at second order. 
Substituting \reff{gla} into Eq.~(\ref{eq1}), all terms at
$\CO(\eps)$ vanish, while at $\CO(\eps^2)$ we obtain equations for
$A_0$ and $A_2$. At $\CO(\eps^3)$ a solvability condition required to
avoid secular growth of $A$ gives the equation 
\huga{\label{gl}
	A_T=c_1(b-b_j) A+c_3|A|^2A, 
} 
where we replaced $\eps$ by
$\sqrt{|b-b_j|}$, and $c_1\equiv\pa_b \lam_1(b_j,k_j)$ and $c_3$ are
the Landau coefficients.  The (analytical) computation of $c_1$ and
$c_3$ is in principle straightforward but can be rather cumbersome,
and we restrict the discussion that follows to their numerical
evaluation using the {\tt ampsys} tool of \pdep\ \cite{ampsys}, and
list their values in Table \ref{Ttab}.  To approximate the Turing
branches we may look for steady (and without loss of generality real)
solutions of \reff{gl}, given by $A=0$ and 
\huga{\label{glsol} A=\sqrt{-c_1(b-b_j)/c_3}, 
}
with $b>b_j$ if $c_1/c_3<0$ (the subcritical regime in our case) and
$b<b_j$ if $c_1/c_3>0$ (the supercritical regime). For illustration, in
Fig.~\ref{f5}(c) we compare $|u_1-u^*(b)|$ obtained from numerical
continuation of the Turing branches (panel (a)) with their first order
Landau approximations, i.e., with $A$ from \reff{glsol}; 
this means a comparison with $u^*(b_j)+2A\cos(k_j x)-u^*(b)$, and since 
$\sup_x \cos(kx)=1$, the dotted branches are given by $2A-(b-b_j)$.

The key reason for the unexpected snaking revealed in panel (a) 
is the strong sensitivity of $c_3$ to $k$ and the sign change 
of $k\mapsto c_3(k)$ near $k\approx 0.42$ 
(see the green branches T5 and T8 in panel (b)), together with 
the nonmonotonic order of the primary BPs, with branches of smaller
wave number (solution B in (a)) interspersed with branches of higher
wave number (solution A in (a)). 
This effectively yields bistability in the
present system despite the fact that the primary branch T1 is supercritical.
Figure~\ref{f5}(a) shows that the low wave number branches T5 and T8
are strongly {\em subcritical} and hence extend substantially past the primary
T1 Turing bifurcation.  Of these the T8 branch (green) extends to the
largest values of the parameter $b$ and yields the widest interval of
bistability with stable $s_4$ and this is indeed the branch on which
the orange snaking branch terminates.

\begin{table*}
	{ 
		\btab{c|ccccccccc}{
			Name&T1&T2&T3&T4&T5&T6&T7&T8&T9\\
			color&blue&red&mag&blue&green&mag&blue&green&mag\\
			$b|_{BP}$&0.9846&0.98452&0.98451&0.98431&0.98414&0.98399&0.98358&0.98337&0.9831\\
			$n$&12&12.5&11.5&13&11&13.5&14&10.5&14.5\\
			$k$&$k_c$&0.46927&0.43172&0.48804&0.41295&0.50681&0.52558&0.39418& 0.54435\\
			$c_1$&-1.0078&-1.0040&-1.0124&-1.0006&-1.0181&-0.9973&-0.9946&-1.02534&-0.9917\\
			$c_3$&-4.1310&-5.8033&-1.7144&-6.9483&1.8053&-7.7145&-8.2141&7.2065&-8.5241	
	}}
	\caption{{\sm Data for the first 9 Turing branches from Fig.~\ref{f5} for 
			$(\del,\alpha,\beta,a,D){=}(0.16,2,0.5,1,0.029)$, $k_c{=}0.4505$, 
			showing the branching points $b|_{BP}$ for solutions with wave numbers
			$k=k_cn/12$, where $n$ is the number of wavelengths on a domain of length
			$\ell=24\pi/k_c$ together with the corresponding Landau coefficients 
			$c_1$ and $c_3$. }
		\label{Ttab}}
\end{table*}

\section{Zipping up a TH snake in the Gilad--Meron model}
In this section we demonstrate that the zipping process in 
Figs.~\ref{f1}-\ref{f3b2} also occurs in another RD system and
can thus be considered generic. We consider the so-called simplified Gilad-Meron
(sGM) model for vegetation patterns in drylands with sandy soil where surface 
water flow is not a major factor. The model considers 
the biomass of the above-ground vegetation,
denoted by $b(x,t)$, and the soil-water content, represented by
$w(x,t)$, and in dimensionless form reads\cite{HPSM04}
\begin{subequations}\label{mode1}
	\begin{align}
	  \partial_t b&= b w \left(1+\eta b \right)^{2}\left(1-b\right)-b+D\partial_x^2 b,  \label{veg1} \\
	\partial_t w&=p-\displaystyle{\frac{n w}{1+\rho b}}-\alpha b w \left(1+\eta b\right)^{2}+\partial_x^2 w, \label{Veg2} \quad
	\end{align}	
\end{subequations} 
where $\eta$ is the root-to-shoot ratio of the plants, $p$ is the precipitation rate,
$\alpha$ is the water uptake efficiency, $n$ is the
evaporation rate, and $\rho$ describes the reduction  
of evaporation resulting from the presence of biomass. 
A linear and nonlinear stability analysis for \reff{mode1} is available in Ref.\cite{ASP23}, together with 
bifurcation diagrams showing various LS branches and the secondary Hopf bifurcations 
from LS leading to short segments of mixed TH states (Ref.\cite{ASP23}, Fig.~9) 
as in Fig.~\ref{f1} for the system \reff{eq1}. 

\begin{figure*}
\btab{lll}{{\sm (a)}&{\sm (b)}&{\sm (d)}\\[-0mm]
\ig[width=60mm]{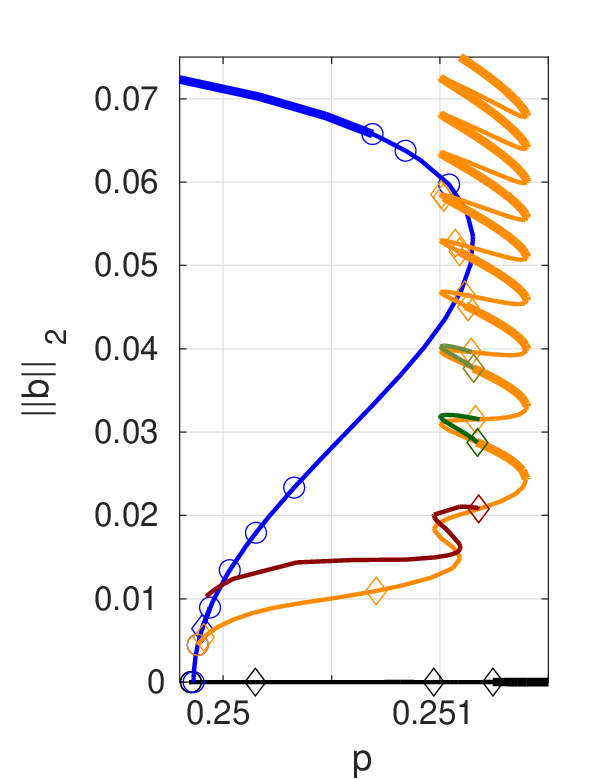}&
\raisebox{35mm}{\btab{l}{\ig[width=57mm]{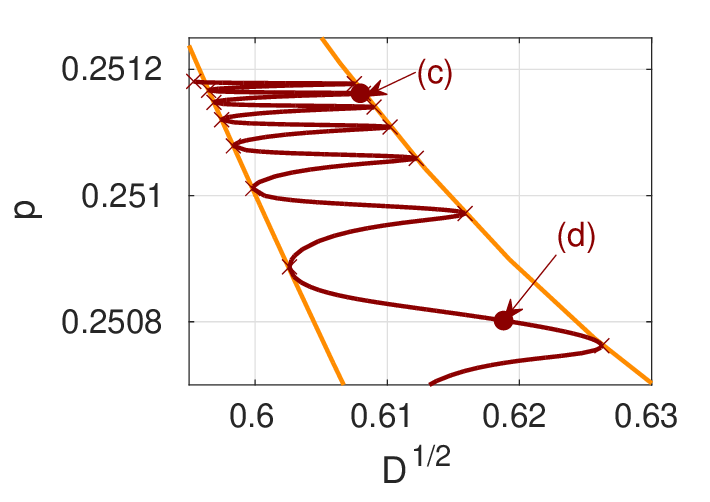}\\[-4mm]
{\sm (c)}\\
\ig[width=60mm]{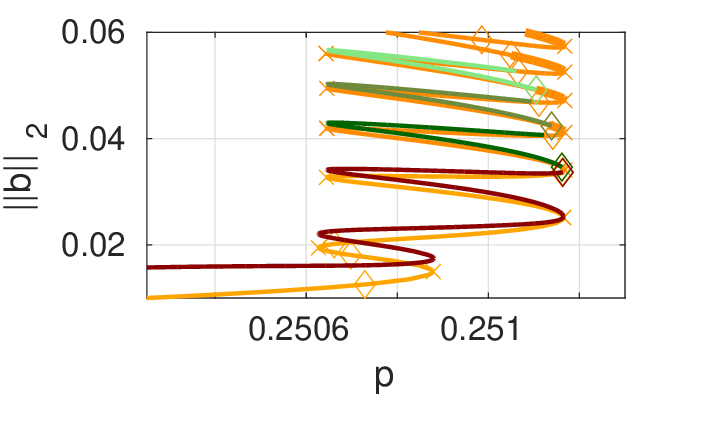}}}
&\ig[width=50mm]{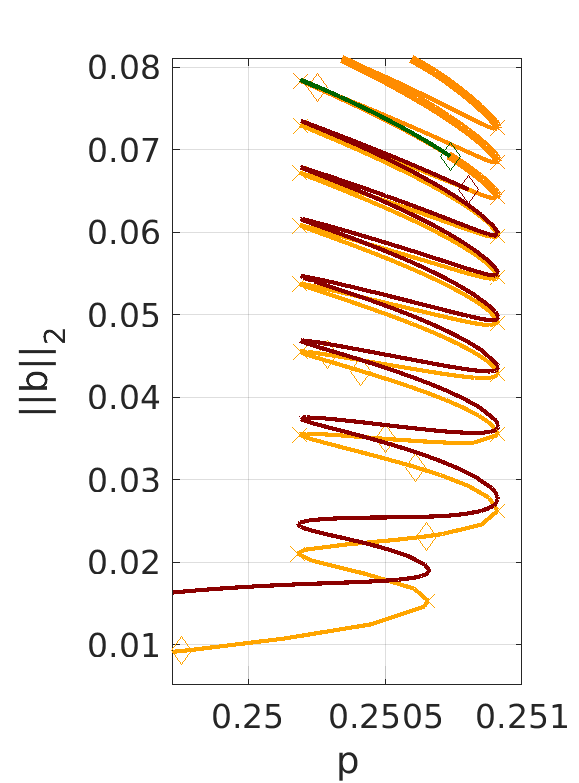}
}
\vs{-4mm}
\caption{{Zipping up the TH snake for the sGM model \reff{mode1} with $\ell{=}60$ and 
    $(n,\eta,\rho,\alpha){=}(0.3579, 3.5, 0.6, 1, 0.5)$. (a) Bifurcation diagram for $D=0.36$ as a
    function of the parameter $p$ showing the primary Turing branch (blue), a steady LS snake
    (orange), and two short TH segments (green). (b) HPC of HP4 (brown), and FPC of FP5
    and FP6 (orange) in the $(D,p)$ plane for comparison with Fig.~\ref{f3b}(b).
    (c) Bifurcation diagram at location (c) in (b), $D=0.37$; (the former) HP4 and HP5 have
    anniliated, and the brown TH branch now starts at HP6. (d) Bifurcation diagram at location
    (d) in (b), $D=0.383$.}\label{snakGM}}
\end{figure*}

Figure \ref{snakGM} shows that we can zip up these segments into a 
TH snake, exactly as in Figs.~\ref{f1}-\ref{f3b2} for the system \reff{eq1}. 
Panel (a) displays the starting situation on a domain slightly smaller ($\ell=60$) than
that in Ref.\cite{ASP23}, Fig.~9. For the ``norm'' $\|b\|_2$ we again use 
\reff{n1},\reff{n2}, where the steady state $b^*$ of uniform vegetation is 
computed numerically at the given parameter values. The first three bifurcations 
from the black branch (at $\|b\|_2=0$) are of Hopf type, and at $p=p_T\approx 0.2499$ 
we find the first Turing branch (T1, blue), 
with a secondary bifurcation to an LS snaking branch at small amplitude (S1, orange). 
The first two HPs on this branch are associated to $k\ne 0$ Hopf 
bifurcations of the background while the third belongs to a $k=0$ Hopf bifurcation 
of the background and the associated brown branch extends down to $p<0.2499$.
The remaining HPs on the LS snake S1 are pairwise connected by short segments, two of which 
are shown in green. Panel (b), similar to Fig.~\ref{f3b}(b), illustrates the HPC (brown)
of HP4 in the parameter $D$ together with the FPC of FP5 and FP6. The HP 
moves through folds given by the FPC, i.e., at these folds two HPs 
collide and annihilate, and the brown branch shows the HP crossing from 
one HP in (a) to the next. In panel (c), HP4 and HP5 have annihilated,
and the brown branch now starts at the (former) HP6, positioned farther
up the orange snake. In panel (d), all the former short green segments
(except for the topmost) have reconnected into a long snaking TH branch.
In contrast to Fig.~\ref{f3b2}, however, this time both the green MM
segments, and the zipped up TH branch (brown) span the same $p$ range as the
LS snake (orange).

\section{Discussion}\label{dsec}
We have revisited the Hopf-Turing interaction that arises
in a number of two-species reaction-diffusion systems.\cite{dewit96,MWBS97,TMBMV13,p2pDMV}
In particular, in Ref.\cite{dewit96} the authors considered the Brusselator model
with supercritical Hopf and Turing branches in a regime with bistability
between the two, finding a large multiplicity of stable Turing states embedded
in an oscillating background obtained via DNS. Since no continuation was
performed the observed states were not linked to homoclinic snaking, a notion
that was developed only subsequently.\cite{CW99} The present work differs in
that we consider a case in which the Turing bifurcation is subcritical. This case,
already considered in Ref.\cite{p2pDMV}, is more interesting since it admits a variety
of steady LS embedded in a homogeneous background. When the Turing bifurcation is
preceded by a (supercritical) Hopf bifurcation, the homogeneous background is no
longer time-independent but begins to oscillate. States of this type were computed
in Ref.\cite{p2pDMV} but their bifurcation behavior was not studied in any detail. The
present work seeks to develop further understanding of these time-periodic LS,
focusing on the model \reff{eq1}. This model, like the Brusselator model studied
in Ref.\cite{dewit96}, also leads to bistability between the Hopf and Turing states,
but this time there is more: the Turing bifurcation is associated with a snaking
branch of steady LS embedded in a trivial state, as also found in Ref.\cite{p2pDMV}.
We have seen that in the bistability region these LS undergo a series of
Hopf bifurcations inherited from the primary Hopf bifurcation that lead to LS
embedded in an oscillating background. We found that for some parameter values
these oscillatory states extend between pairs of Hopf bifurcations on the snaking
LS branch, one on either side of every right fold. Remarkably, we found that by
varying a second parameter we could progressively move a secondary Hopf bifurcation to
a small amplitude mixed mode up the LS branch, leading to repeated reconnection between
the mixed mode and the disconnected oscillatory states on the LS branch. We described
the net effect of these reconnections as {\it zipping up} of these disconnected branches
into a snaking branch of oscillatory states and demonstrated that similar zipping up
occurs in other RD systems such as the simplified Gilad-Meron model. We believe our study
presents perhaps the best example of this behavior in a two-species RD system, and the most
compelling example of snaking of spatially localized time-periodic states.
Similar zipping transitions arise in forced snaking\cite{ponedel16} and are
studied in detail in Ref.\cite{parra23}, albeit for steady LS only. See also
Refs.\cite{verschueren21,tumelty23}. 

Once a TH snake is established as in Fig.~\ref{f3b2} for \reff{eq1} 
by varying $b$, or in Fig.~\ref{snakGM} for \reff{mode1} by varying $p$,
we can choose a point from such branches and find snaking behavior on
continuing that point in some other parameter. Thus the primary parameters
$b$ in \reff{eq1} and $p$ in \reff{mode1} are essentially just convenient 
standard choices. However, in both our models \reff{eq1} and \reff{mode1}, 
we need to use the diffusion constant $D$ for the zipping up of the TH snake,
i.e., to the best of our knowledge a similar zipping up is not possible using
two--parameter continuation that excludes $D$; see Fig.~\ref{f3c} for an example
using $(b,\del)$ in \reff{eq1}. Technically, the explanation for the zipping up
is given in Remark \ref{reconrem}a). This remark applies to \reff{mode1} as well,
but why this occurs in both models under Hopf--point continuation in $D$ but not
other parameters warrants further investigation.

Surprisingly, we also found homoclinic snaking in the regime where the Turing
bifurcation of the trivial state comes in first and is supercritical, a consequence
of a strongly subcritical secondary bifurcation to spatially modulated states
that extend well into the regime of stable trivial states. The resulting LS are
then embedded in the trivial state but when this state becomes Turing-unstable
they connect to a small amplitude Turing state in the background. These LS were
found to snake (and hence acquire stability) owing to bistability with strongly
subcritical sideband Turing states. We traced this behavior to the nonmonotonic
order of the primary bifurcation points with respect to the wave number $k$ together
with the fact that the key Landau coefficient $c_3$ depends sensitively on $k$ and 
changes sign near $k_c$. 
It must be emphasized, however, that whenever the LS are embedded in a periodic
background the imposed domain size plays a significant role. This is in addition
to its role in determining the order of the primary BPs to the various sidebranches
associated with the Turing instability.

Altogether, our results add to the variety, multistability, and competition 
between patterned or spatially localized steady states and time--periodic states
of two--component RD systems. For instance, Fig.~\ref{fdns} illustrates that in
the vicinity of a TH snake very different long--time behavior may be found, depending
both on the details of the initial condition and the parameter value. Moreover, the
zipping up of the even TH snake (and a similar zipping up of the odd TH snake -- not
shown) generates a pair of intertwined snaking branches of time-periodic states and
hence a large multiplicity of coexisting time-periodic states of either parity, in
complete analogy with the snaking of steady LS, thereby explaining why LS embedded
in an oscillatory background may occur in systems with arbitrarily large spatial extent.

There is a natural extension of the present work, and that is to three-species RD
systems. These systems may exhibit, in appropriate parameter regimes, a wave instability,
i.e., a Hopf bifurcation with a finite wave number. This instability may develop into
standing or traveling waves, depending on parameters.\cite{knobloch86} We anticipate
therefore that in such systems we may be able to compute LS embedded in a background
of standing or traveling waves, depending on parameters, as already found for Marangoni
driven convection.\cite{assemat08} In two-species systems such a bifurcation is never
the first instability for which $k_H=0$ but we have seen that $k_H\ne 0$ states can
set in in subsequent primary Hopf bifurcations, and that this bifurcation can likewise
be inherited by the secondary LS, resulting in spatially localized time-periodic structures. It
is of interest to determine whether these states also snake. We mention that numerically
stable steady LS are found even when the background trivial state is unstable to
traveling waves, provided it is only convectively unstable.\cite{batiste06,watanabe}
Thus the three-species case is expected to be considerably richer than the system \reff{eq1}
studied here.

\vs{3mm}
\noi
{\bf Acknowledgement}. The work of EK was supported in part by the National Science Foundation
under Grant No. DMS-1908891.

\bibliographystyle{plain}
\bibliography{Ref}


\end{document}